\documentclass[12pt]{amsart}
\usepackage{amssymb,amsmath}
\usepackage{graphicx}
\usepackage{colordvi}

\numberwithin{equation}{section}

\newtheorem{theorem}[equation]{Theorem}

\newtheorem{proposition}[equation]{Proposition}
\newtheorem{corollary}[equation]{Corollary}

\theoremstyle{definition}
\newtheorem{definition}[equation]{Definition}

\newtheorem{example}[equation]{Example}

\newtheorem{remark}[equation]{Remark}

\newcommand{\DeCo}[1]{\BlueViolet{#1}}
\newcommand{\Sim}{\includegraphics[height=7pt]{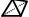}}
\newcommand{\bSim}{\includegraphics[height=7pt]{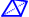}}
\newcommand{\tri}{\includegraphics{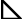}}
\newcommand{\btri}{\includegraphics{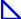}}

\newcommand{\C}{{\mathbb C}}
\newcommand{\R}{{\mathbb R}}
\newcommand{\Z}{{\mathbb Z}}

\newcommand{\bb}{{\bf b}}

\begin{document}

\title[Geometrical aspects of control points]{Some geometrical aspects of control points\\
  for toric patches} 

\author[G.~Craciun]{Gheorghe Craciun}
\address{Department of Mathematics\\
         and Department of Biomolecular Chemistry\\
         University of Wisconsin\\
         Madison\\
         WI \ 53706\\
         USA}
\email{craciun@math.wisc.edu}
\urladdr{www.math.wisc.edu/\~{}craciun}
\author[L.~Garc{\'\i}a-Puente]{Luis David Garc{\'\i}a-Puente}
\address{Department of Mathematics and Statistics\\
         Sam Houston State University\\
         Huntsville\\
         TX \ 77341\\
         USA}
\email{lgarcia@shsu.edu}
\urladdr{www.shsu.edu/\~{}ldg005}
\author[F.~Sottile]{Frank Sottile}
\address{Department of Mathematics\\
         Texas A\&M University\\
         College Station\\
         TX \ 77843\\
         USA}
\email{sottile@math.tamu.edu}
\urladdr{www.math.tamu.edu/\~{}sottile}
\thanks{This material is based in part upon work
   supported by the Texas Advanced Research
   Program under Grant No. 010366-0054-2007 and  NSF grant DMS-070105.} 


\begin{abstract}
  We use ideas from algebraic geometry and dynamical
  systems to explain some ways that control points influence the shape of a B\'ezier 
  curve or patch. 
  In particular, we establish a generalization of Birch's Theorem and use it to deduce
  sufficient conditions on the control points for a patch to be injective.
  We also explain a way that the control points influence the shape via degenerations 
  to regular control polytopes.
  The natural objects of this investigation are irrational patches, which are 
  a generalization of Krasauskas's toric patches, and include B\'ezier and tensor product
  patches as important special cases. 
\end{abstract}
\maketitle              

%
\section*{Introduction}
%

The control points and weights of a B\'ezier curve, B\'ezier patch, or tensor-product
patch govern many aspects of the curve or surface.
For example, they provide an intuitive means to control its shape.
Through de Castlejau's algorithm, they enable the computation of 
the curve or surface patch.
Finer aspects of the patch, particularly continuity and smoothness at the boundary 
of two patches are determined by the control points and weights.
Global properties, such as the location of a patch in space due to the convex hull
property, also depend upon the control points.
When the control points are in a particular convex position, then the patch is
convex~\cite{CGS-Da91}. 

We apply methods from algebraic geometry, specifically toric geometry, to explain
how some further global properties of a patch are governed by the control points.
We first investigate the self-intersection, or injectivity of a patch. 
We give a simple and easy-to-verify condition on a set of control points which implies
that the resulting patch has no self-intersection, for any choice of weights.
For 3-dimensional patches as used for solid modeling, injectivity is equivalent
to the patch properly parameterizing the given solid.
This uses Craciun and Feinberg's injectivity theorem~\cite{CGS-CF05}
from the theory of chemical reaction networks, which may be seen as a generalization of
Birch's Theorem from algebraic statistics. 

A second global property that we investigate is how the shape of the patch is related to
the shape of a \DeCo{{\sl control polytope}}.
This is a piecewise linear triangulated surface whose vertices are the control points.
It is \DeCo{{\sl regular}} if the underlying triangulation comes from a regular
triangulation of the domain polytope of the patch.
We show that regular control polytopes are the limits of patches as the weights
undergo a toric deformation corresponding to the underlying regular triangulation,
and that non-regular control polytopes can never be such a limit.
This gives a precise meaning to the notion that the shape of the control net governs the
shape of the patch.

This line of inquiry is pursued in terms of Krasauskas's toric
patches~\cite{CGS-Kr02}, as it relies upon the structure of toric varieties from
algebraic geometry.
The correct level of generality is however that of \DeCo{{\sl irrational (toric) patches}},
which are analytic subvarieties of the simplex (realized as a compactified positive orthant) 
that are parameterized by monomials \DeCo{$x^\alpha$}, where \DeCo{$x$} is a vector of
positive numbers and the exponent vector \DeCo{$\alpha$} has real-number coordinates.
(This is a usual toric variety when $\alpha$ has integer coordinates.)
While irrational patches may seem exotic for modeling, they occur naturally in statistics
as discrete exponential families~\cite{CGS-Br86} and their blending functions may be computed
using iterative proportional fitting (IPF)~\cite{CGS-DR72}, a popular numerical algorithm from 
statistics.
Furthermore, these blending functions have linear precision.
For toric patches, this was observed in~\cite{CGS-So03} and developed in~\cite{CGS-GS}, and the
analysis there carries over to irrational patches.
While we work in this generality, our primary intent (and the main application) is to shed
light on properties of B\'ezier curves, surfaces, and 3-dimensional patches.

We recall the standard definition of a mapping via control points and blending functions,
and then the definitions of toric B\'ezier patches in Sect.~\ref{CGS-S:defs}.
There, we also illustrate some of our results on examples of B\'ezier curves.
In Sect.~\ref{CGS-S:irrational}, we introduce irrational toric patches, recalling the
geometric formulation of a toric patch and the use of iterative proportional fitting to
compute these patches, explaining how these notions from~\cite{CGS-GS} for toric patches
extend to irrational patches. 
The next two sections contain our main results.
We study injectivity of patches in Sect.~\ref{CGS-S:inj}, and discuss degenerations to
control polytopes in Sect.~\ref{CGS-S:degeneration}.
Appendices A and B contain technical proofs of some theorems.

%
\section{Toric B\'ezier Patches}\label{CGS-S:defs}
%

We interpret the standard definition of a mapping via control
points and blending functions (see for example \cite[\S 2]{CGS-KK00}) 
in a general form convenient for our discussion.
All functions here are smooth ($C^\infty$) where defined and
real-valued.
Let \DeCo{$\R_>$} be the set of strictly positive real numbers and \DeCo{$\R_\geq$} the 
set of  non-negative real numbers.
We will use the following typographic conventions throughout.
Vector constants (control points, indexing exponents, and standard basis vectors) will be 
typeset in {\bf  boldface}, while vector variables will be typeset in standard math italics.

Let ${\mathcal A}$ be a finite set of points that affinely span $\R^d$, which we shall use as
geometrically meaningful indices. 
A control point scheme for parametric patches, or 
\DeCo{({\sl parametric}) {\sl patch}}, is a collection
$\beta=\{\beta_{\bf a} \mid {\bf a}\in{\mathcal A}\}$ of non-negative functions, called
\DeCo{{\sl blending functions}}.
The common domain of the blending functions is the convex hull $\Delta$ of
${\mathcal A}$, which we call the \DeCo{{\sl domain polytope}}.
We also assume that the blending functions do not vanish simultaneously 
at any point of $\Delta$, so that there are no basepoints.

Lists $\DeCo{{\mathcal B}}:=\{{\bf b}_{\bf a} \mid {\bf a}\in{\mathcal A}\}\subset\R^n$ of 
\DeCo{{\sl control points}}
and positive \DeCo{{\sl weights}}
$\DeCo{w}:=\{w_{\bf a}\in\R_>\mid{\bf a}\in{\mathcal A}\}\in\DeCo{\R^{\mathcal A}_>}$   
together give a map  $F \colon\Delta\to\R^n$ defined by
 \begin{equation}\label{CGS-Eq:P_Patch}
   F(x)\ :=\ 
   \frac{\sum_{{\bf a}\in{\mathcal A}} w_{\bf a}\beta_{\bf a}(x)\,{\bf b}_{\bf a} }
        {\sum_{{\bf a}\in{\mathcal A}} w_{\bf a}\beta_{\bf a}(x)}\enspace .
 \end{equation}
The denominator in~\eqref{CGS-Eq:P_Patch} is positive on $\Delta$ and so
the map $F$ is well-defined. 

\begin{remark}
 We will refer to both $\beta_{\bf a}(x)$ and $w_{\bf a} \beta_{\bf a}(x)$ as the
 blending functions of a patch.
 This generality of separating the weights from the blending functions
 will be used in Sect.~\ref{CGS-S:degeneration} when we investigate the effect of
 systematically varying the weights of a patch while keeping the control points and
 blending functions constant. 
\end{remark}

The control points and weights affect the shape of the patch which is the image of the map
$F$~\eqref{CGS-Eq:P_Patch}. 
For example, the \DeCo{{\sl convex hull property}} asserts that 
the image $F(\Delta)$ of the patch lies in the convex hull of the control points. 
To see this, note that if we set
\[
   \DeCo{\overline{\beta}_{\bf a}(x)}\ :=\ 
    \frac{w_{\bf a}\beta_{\bf a}(x)}{\sum_{{\bf a}\in{\mathcal A}} w_{\bf a}\beta_{\bf a}(x)}\enspace,
\]
then $\overline{\beta}_{\bf a}(x)\geq 0$ and 
$1=\sum_{{\bf a}\in{\mathcal A}} \overline{\beta}_{\bf a}(x)$.
Then formula~\eqref{CGS-Eq:P_Patch} becomes
\[
   F(x)\ =\ 
   \sum_{{\bf a}\in{\mathcal A}} \overline{\beta}_{\bf a}(x) \,{\bf b}_{\bf a} \enspace ,
\]
so that $F(x)$ is a convex combination of the control points and therefore lies in their
convex hull.
In fact, if there is a point $x\in\Delta$ at which no blending function vanishes, then 
{\sl any} point in the interior of the convex hull of the control
points is the image $F(x)$ of some patch for some choice of weights.
In this way, the convex hull property is the strongest general statement that can be
made about the location of a patch.

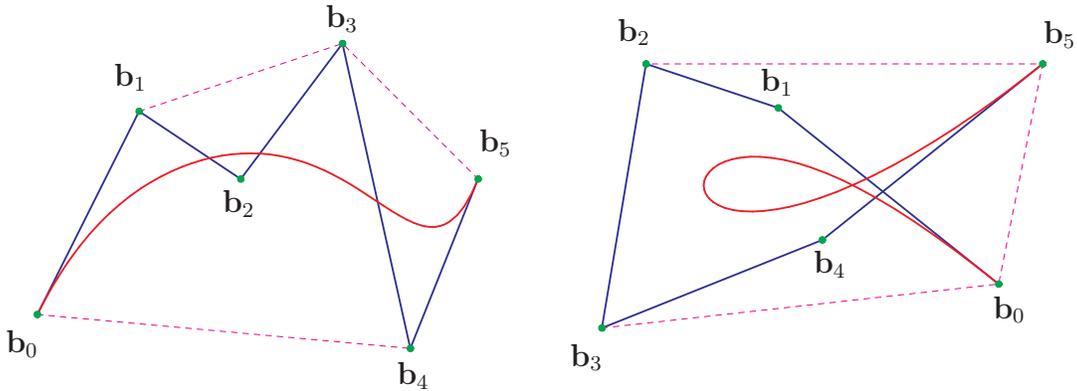
\begin{figure}[htb]
 \[
  \begin{picture}(190,150)(-9,-10)
     \put(0,0){\includegraphics{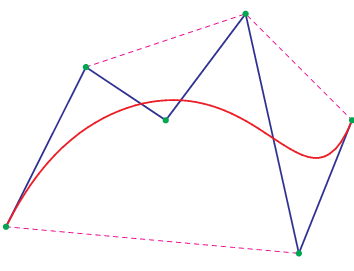}}
     \put(170,83){$\bb_5$}
     \put(139,-5){$\bb_4$}
     \put(112,130){$\bb_3$}
     \put(73,60){$\bb_2$}
     \put(32,108){$\bb_1$}
     \put(-9,7){$\bb_0$}
   \end{picture}
    \qquad
   \begin{picture}(190,150)(-9,-10)
     \put(0,0){\includegraphics{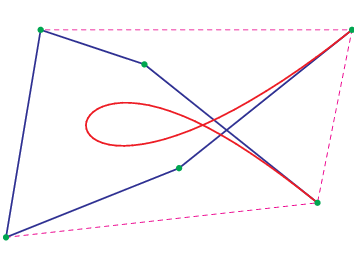}}
     \put(170,124){$\bb_5$}
     \put( 83, 37){$\bb_4$}
     \put( -9,  2){$\bb_3$}
     \put(  9,125){$\bb_2$}
     \put( 63,106){$\bb_1$}
     \put(151, 20){$\bb_0$}
  \end{picture}
 \]
 \caption{Quintic B\'ezier curves.}\label{CGS-F:one}
\end{figure}

Another well-known manifestation of control points is the relation of a B\'ezier curve to
its control polygon.
Fix a positive integer $m$ and let 
${\mathcal A}:= \{\frac{i}{m} \mid i=0,\dotsc,m\}$ so that $\Delta$ is the unit interval. 
The blending functions of a B\'ezier curve are the Bernstein polynomials,
\[
   \DeCo{\beta_i(x)}\ (=\ \beta_{\frac{i}{m}}(x))\ :=\ 
   \tbinom{m}{i} x^i(1-x)^{m-i}\enspace.
\]
The \DeCo{{\sl control polygon}} of a B\'ezier curve with control points 
${\bf b}_0,{\bf b}_1,\dotsc,{\bf b}_m$ is
the union of the line segments 
$\overline{{\bf b}_0,{\bf b}_1},\, \overline{{\bf b}_1,{\bf b}_2},\,\dotsc,\,
  \overline{{\bf b}_{m-1},{\bf b}_m}$
between consecutive control points.
Figure~\ref{CGS-F:one} displays two quintic plane B\'ezier curves with their control
polygons (solid lines).  
The convex hulls of the control points are indicated by the dashed lines.
The first curve has no points of self-intersection, while the
second curve has one point of self-intersection.
While this self-intersection may be removed by varying the weights attached to the control
points, by Theorem~\ref{CGS-Th:inj} it is impossible to find weights so
that a curve with the first set of control points has a point of self-intersection.

We will also show that the control polygon may be approximated by
a B\'ezier curve.
We state a simplified version of Theorem~\ref{CGS-Th:Bez_curve_approx} from
Sect.~\ref{CGS-S:degeneration}. \medskip

\noindent{\bf Theorem.} {\it
  Given control points in $\R^n$ for a B\'ezier curve and some number $\epsilon>0$,
  there 
  is a choice of weights so that the image $F[0,1]$ of the B\'ezier curve lies within
  a distance $\epsilon$ of the control polygon.}\medskip

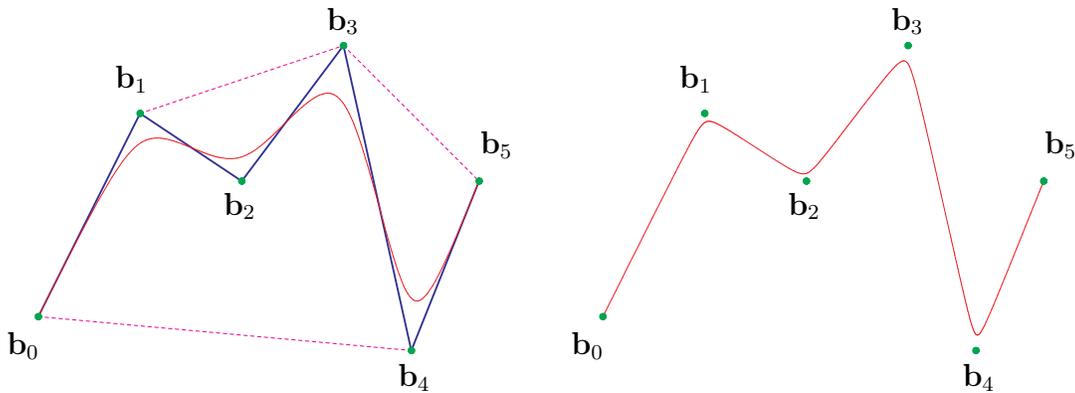
\begin{figure}[htb]
 \[
  \begin{picture}(190,150)(-9,-10)
     \put(0,0){\includegraphics{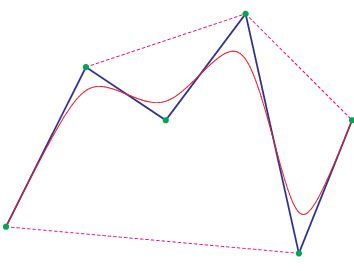}}
     \put(170,83){$\bb_5$}
     \put(139,-5){$\bb_4$}
     \put(112,130){$\bb_3$}
     \put(73,60){$\bb_2$}
     \put(32,108){$\bb_1$}
     \put(-9,7){$\bb_0$}
  \end{picture}
  \qquad
  \begin{picture}(190,150)(-9,-10)
     \put(0,0){\includegraphics{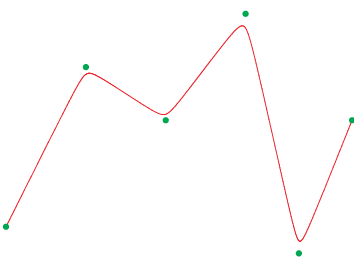}}
     \put(170,83){$\bb_5$}
     \put(139,-5){$\bb_4$}
     \put(112,130){$\bb_3$}
     \put(73,60){$\bb_2$}
     \put(32,108){$\bb_1$}
     \put(-9,7){$\bb_0$}
  \end{picture}
 \]
 \caption{Degenerating quintics.}\label{CGS-F:Degenerating}
\end{figure}

In Figure~\ref{CGS-F:Degenerating}, we display one of the quintic curves from Figure~\ref{CGS-F:one}, 
but with weights on ${\bf b}_0$---${\bf b}_5$ of $(1,20^2,20^3,20^3,20^2,1)$ and 
$(1,300^2,300^3,300^3,300^2,1)$, respectively.
The control polygon for the second curve is omitted, as it would obscure the curve.
The first curve lies within a distance $\epsilon=0.13$ of the control polygon and the 
second within a distance $\epsilon=0.02$, if the control polygon has height 1. 

%
%
\subsection{Toric Patches}
 Krasauskas~\cite{CGS-Kr02} introduced toric patches as a generalization of the classical
 B\'ezier and tensor product patches.
 These are based upon toric varieties from algebraic geometry and 
 their shape may be 
 any polytope with integer vertices.
 The articles~\cite{CGS-Cox03,CGS-So03} provide an introduction to toric varieties for 
 geometric modeling.

 A polytope $\Delta$ is defined by its \DeCo{{\sl facet inequalities}}
\[
  \Delta\ =\ \{x\in\R^d \mid 0 \leq h_i(x)\,, i=1,\dotsc,\ell\}\enspace.
\]
 Here, $\Delta$ has $\ell$ facets (faces of maximal dimension) and for each
 $i=1,\dotsc,\ell$,  $h_i(x)={\bf v}_i\cdot x + c_i$ is the linear function defining
 the $i$th facet, where ${\bf v}_i\in\Z^d$ is the (inward oriented) primitive vector normal to 
 the facet and $c_i\in\Z$. 

 For example, if our polytope is the triangle with vertices $(0,0)$, $(m,0)$, and $(0,m)$,  
 \begin{equation}\label{CGS-Eq:d-tri}
   \DeCo{m}\,\btri\ :=\  
    \{(x,y)\ \in\ \R^2 \mid 0\leq x,y, \quad\mbox{and}\quad 0\leq m-(x+y)\}\enspace,
 \end{equation}
 then we have $h_1=x$, $h_2=y$, and $h_3=m-x-y$.
 Here, $m\,\tri$ is the unit triangle $\tri$ with vertices $(0,0)$, $(1,0)$, and $(0,1)$
 scaled by a factor of $m$. 

 Let ${\mathcal A}\subset \Delta\cap\Z^d$ be any subset of the integer points of
 $\Delta$ which includes its vertices.
 For every ${\bf a}\in{\mathcal A}$, Krasauskas defined the \DeCo{{\sl toric B\'ezier function}}
 \begin{equation}\label{CGS-Eq:toric-Bezier}
   \DeCo{\beta_{\bf a}(x)}\ :=\
    h_1(x)^{h_1({\bf a})}h_2(x)^{h_2({\bf a})}\dotsb h_\ell(x)^{h_\ell({\bf a})}\enspace,
 \end{equation}
 which is non-negative on $\Delta$, and the collection of all $\beta_{\bf a}$ has no common
 zeroes on $\Delta$. 
 These are blending functions for the \DeCo{{\sl toric patch}} of shape
 ${\mathcal A}$.
 If we choose weights $w\in\R^{\mathcal A}$ and multiply the
 formula~\eqref{CGS-Eq:toric-Bezier} by $w_{\bf a}$, we obtain blending functions for the toric
 patch of shape $({\mathcal A},w)$.

\begin{example}[B\'ezier triangles]\label{CGS-Ex:BezTri}
 When $\Delta$ is a scaled triangle or a product of such
 triangles,~\eqref{CGS-Eq:toric-Bezier} gives the blending 
 functions of the B\'ezier patch or B\'ezier simploid~\cite{CGS-dRGHM} with the corresponding
 shape.
 To see this for the scaled triangle $m\,\tri$~\eqref{CGS-Eq:d-tri}, note that
 given an integer point ${\bf a}=(i,j)\in m\,\tri$,
 and weight the multinomial coefficient $w_{(i,j)}:=\frac{m!}{i!j!(m-i-j)!}$,
 then the corresponding 
 blending function is
\[
   \beta_{(i,j)}(x,y)\ =\  \tfrac{m!}{i!j!(m-i-j)!} x^i y^j (m-x-y)^{m-i-j}\enspace.
\]
 This is almost the bivariate Bernstein polynomial, which is obtained by 
 substituting $mx$ and $my$ for $x$ and $y$, respectively, and dividing by $m^m$.
 (This has the effect of changing the domain from $m\,\tri$ to the unit triangle $\tri$.)
\end{example}

%
%
\section{Irrational Patches}\label{CGS-S:irrational}

Krasauskas's definition~\eqref{CGS-Eq:toric-Bezier} of toric B\'ezier functions still makes
sense if we relax the requirement that the points ${\mathcal A}\subset\R^d$ have integer
coordinates.
This leads to the notion of an irrational patch (as its blending functions are no longer
rational functions), which provides the level of generality appropriate for our
investigation. 

Let ${\mathcal A}\subset\R^d$ be a finite collection of points and set $\Delta\subset\R^d$ to be
the convex hull of ${\mathcal A}$, which we assume is a full-dimensional polytope.
We may also realize $\Delta$ as an intersection of half-spaces through its facet
inequalities,
 \begin{equation}\label{CGS-Eq:facet}
    \Delta\ =\ \{x\in\R^d \mid h_i(x)\geq 0\ \mbox{for each}\ i=1,\dotsc,\ell\}\enspace.
 \end{equation}
Here, $\Delta$ has $\ell$ facets with the $i$th facet supported by the affine hyperplane
$h_i(x)=0$ where $h_i(x)={\bf v}_i\cdot x+c_i$ with $c_i\in\R$ and ${\bf v}_i$ {\it an} inward
pointing normal vector to the $i$th facet of $\Delta$. 
There is no canonical choice for these data; multiplying a pair 
$({\bf v}_i,c_i)$ by a positive scalar gives another pair defining the same half-space. 

Following Krasauskas, we provisionally define  
\DeCo{{\sl (irrational)  toric B\'ezier functions}} 
$\{\beta_{\bf a}\colon\Delta\to\R_\geq \mid {\bf a}\in{\mathcal A}\}$
by the same formula as~\eqref{CGS-Eq:toric-Bezier},
\[
   \beta_{\bf a}(x)\ :=\
    h_1(x)^{h_1({\bf a})}h_2(x)^{h_2({\bf a})}\dotsb h_\ell(x)^{h_\ell({\bf a})}\enspace.
\]
These are blending functions for the \DeCo{{\sl irrational toric patch}} of shape ${\mathcal A}$.
While these functions do depend upon the choice of data $({\bf v}_i,c_i)$ for the facet
inequalities defining $\Delta$, we will see that the image $F(\Delta)$ of such a patch given by weights and
control points is independent of these choices. 

%
%
\subsection{Geometric Formulation of a Patch}
We follow Sect.~2.2 of~\cite{CGS-GS}, but drop the requirement that our objects are
algebraic.
Let ${\mathcal A}\subset\R^d$ be a finite subset indexing a collection of blending functions
$\{\beta_{\bf a}\colon\Delta\to\R_\geq\mid{\bf a}\in{\mathcal A}\}$, where $\Delta$ is the convex
hull of ${\mathcal A}$. 
Let $\R^{\mathcal A}$ be a real vector space with basis $\{e_{\bf a}\mid{\bf a}\in{\mathcal A}\}$.
Set $\R^{\mathcal A}_\geq\subset\R^{\mathcal A}$ to be the points with non-negative coordinates and let
$\R^{\mathcal A}_>$ be those points with strictly positive coordinates.

For $z=(z_{\bf a} \mid {\bf a}\in{\mathcal A})\in\R^{\mathcal A}_\geq$, set
$\DeCo{\sum z}:=\sum_{{\bf a}\in{\mathcal A}} z_{\bf a}$.
The \DeCo{{\sl ${\mathcal A}$-simplex}}, $\bSim^{\mathcal A}\subset\R^{\mathcal A}_\geq$, is the set
\[
  \bSim^{\mathcal A}\ :=\ 
   \{ z\ \in\ \R^{\mathcal A}_{\geq} :  {\textstyle \sum z}=1\}\enspace.
\]
We introduce homogeneous coordinates for $\Sim^{\mathcal A}$.
If $z\in\R^{\mathcal A}_\geq\setminus\{0\}$, then we set
\[
   [z_{\bf a} \mid  {\bf a}\in{\mathcal A}]\ :=\ \frac{1}{\sum z}(z_{\bf a} \mid {\bf a}\in{\mathcal A})\ \in\ \Sim^{\mathcal A}\enspace.
\]
The blending functions $\{\beta_{\bf a}\colon\Delta\to\R_\geq\mid{\bf a}\in{\mathcal A}\}$ give a 
$C^\infty$ map, 
 \begin{eqnarray*}
    \beta\ \colon\ \Delta&\longrightarrow& \Sim^{\mathcal A}\\
                        x&\longmapsto&[\beta_{\bf a}(x) \mid {\bf a}\in{\mathcal A}]\,.
 \end{eqnarray*}

The reason for this definition is that a mapping
$F\colon\Delta\to\R^n$~\eqref{CGS-Eq:P_Patch} 
given by the blending functions $\beta$, weights  $w$, 
and control points ${\mathcal B}$
factors through the map $\beta\colon\Delta\to\Sim^{\mathcal A}$.
To see this, first note that the weights $w\in\R^{\mathcal A}_>$ act on $\Sim^{\mathcal A}$:
If $z=[z_{\bf a} \mid {\bf a}\in{\mathcal A}]\in \Sim^{\mathcal A}$, then 
 \begin{equation}\label{CGS-Eq:w-map}   
    \DeCo{w.z}\ :=\ [w_{\bf a} z_{\bf a}  \mid {\bf a}\in{\mathcal A}]\enspace.
 \end{equation}
The control points ${\mathcal B}$ define the map $\pi_{\mathcal B}\colon\Sim^{\mathcal A}\to\R^n$
via
 \[
    \pi_{\mathcal B}\ \colon\ z\ =\ (z_{\bf a} \mid {\bf a}\in{\mathcal A})\ \longmapsto\ 
    \sum_{{\bf a}\in{\mathcal A}} z_{\bf a} {\bf b}_{\bf a}\enspace .
 \]
Then the mapping $F$~\eqref{CGS-Eq:P_Patch} is simply the composition
 \begin{equation}\label{CGS-Eq:composition}
   \Delta\ \xrightarrow{\;\beta\;}\ \Sim^{\mathcal A}\ \xrightarrow{\;w.\;}\ 
    \Sim^{\mathcal A}\ \xrightarrow{\;\pi_{\mathcal B}\;}\ \R^n\enspace.
 \end{equation}

In this way, we see that the image $\beta(\Delta)\subset\Sim^{\mathcal A}$ of $\Delta$ under the map
$\beta$ determines the shape of the patch $F(\Delta)$~\eqref{CGS-Eq:P_Patch}.
Internal structures of the patch, such as the mapping of texture, are determined by how
$\beta$ maps $\Delta$ to $\beta(\Delta)$. For example, precomposing $\beta$ with any
homeomorphism of $\Delta$ gives blending functions with the same image in $\Sim^{\mathcal A}$, but
with a different internal structure.

For an irrational toric patch of shape ${\mathcal A}$ with blending
functions~\eqref{CGS-Eq:toric-Bezier}, the image $\beta(\Delta)\subset\Sim^{\mathcal A}$ is independent of
the choice of normal vectors.
For this, we first define the map $\DeCo{\varphi_{\mathcal A}}\colon\R^d_>\to \Sim^{\mathcal A}$ by
 \begin{equation}\label{CGS-Eq:monom_parametrize}
   \varphi_{\mathcal A}\ \colon\ 
     (x_1,\dotsc,x_d)\ \longmapsto\ [ x^{\bf a}\ :\ {\bf a}\in{\mathcal A}]\enspace.
 \end{equation}
Let \DeCo{$X_{\mathcal A}$} be the closure of the image of the map $\varphi_{\mathcal A}$.
When ${\mathcal A}\subset\Z^d$, this is the positive part~\cite[\S 4]{CGS-Fu93} of the toric variety
parameterized by the monomials of ${\mathcal A}$.
When ${\mathcal A}$ is not integral, we call $X_{\mathcal A}$ the 
\DeCo{{\sl (irrational) toric variety}} parameterized by monomials in ${\mathcal A}$.

In Appendix~\ref{CGS-A:B} we prove the following theorem.

\begin{theorem}\label{CGS-Th:X_calA}
  Suppose that ${\mathcal A}\subset\R^d$ is a finite set of points with convex hull $\Delta$.
  Let $\beta=\{\beta_{\bf a}\mid {\bf a}\in{\mathcal A}\}$ be a collection of irrational toric B\'ezier
  functions for ${\mathcal A}$. 
  Then $\beta(\Delta)=X_{\mathcal A}$,  the closure of the image of $\varphi_{\mathcal A}$.
\end{theorem}

We prove this by showing that the restriction of the map $\beta$ to the interior
$\Delta^\circ$ of $\Delta$ factors through the map $\varphi_{\mathcal A}$.

By Theorem~\ref{CGS-Th:X_calA}, the image of the irrational toric blending
functions for ${\mathcal A}$ depends upon ${\mathcal A}$ and not upon the choice of toric blending
functions for ${\mathcal A}$.
Thus the shape of the corresponding patch $F(\Delta)$~\eqref{CGS-Eq:P_Patch} depends only upon
${\mathcal A}$, the weights $w$, and the control points ${\mathcal B}$.
However, the actual parameterization of $X_{\mathcal A}$ by $\Delta$, and hence of $F(\Delta)$
does depend upon the choice of toric blending functions for ${\mathcal A}$.

To ensure that the patch $F(\Delta)$ has shape reflecting that of $\Delta$, we require
that the map $\beta\colon \Delta \to X_{\mathcal A}$ be injective.
This also guarantees that the patch $F(\Delta)$ is typically an immersion.
In the context of irrational toric patches, this injectivity is guaranteed by Birch's
Theorem from algebraic statistics.  
For a standard reference, see~\cite[p.~168]{CGS-AGR90}.
When ${\mathcal A}\subset\Z^d$, Birch's Theorem follows from general results on the
moment map in symplectic geometry~\cite[\S 4.2 and Notes to Chapter 4, p.~140]{CGS-Fu93}.

\begin{theorem}[Birch's Theorem]\label{CGS-Th:Birch}
  Suppose ${\mathcal A}\subset\R^d$ is finite and let $\beta$ be a collection of toric B\'ezier
  functions. 
  If we choose control points ${\mathcal B}$ to be the corresponding points of ${\mathcal A}$, 
  $\{{\bf b}_{\bf a}={\bf a}\mid{\bf a}\in{\mathcal A}\}$, then the composition
\[
  \Delta\ \xrightarrow{\ \beta\ }\ X_{\mathcal A}\ 
          \xrightarrow{\ \pi_{\mathcal B}\ }\ \R^d
\]
 is a homeomorphism onto $\Delta$.
\end{theorem}

By Birch's Theorem and Theorem~\ref{CGS-Th:X_calA}, any two sets $\beta,\beta'$ of toric
B\'ezier functions of shape ${\mathcal A}$ differ only by a homeomorphism 
$h\colon\Delta\xrightarrow{\ \sim\ }\Delta$ of the polytope $\Delta$, so that
$\beta'=\beta\circ h$.
In fact $h$ restricts to a homeomorphism on all faces of $\Delta$.
As we are concerned with the shape of a patch and not its internal structure, we follow
Krasauskas' lead and make the following definition.

\begin{definition}
  A(n irrational) toric patch of shape ${\mathcal A}$ is any set of blending functions
  $\beta:=\{\beta_{\bf a}\colon\Delta\to\R_\geq\mid {\bf a}\in{\mathcal A}\}$ such that the map
  $\beta\colon\Delta\to X_{\mathcal A}$ is a homeomorphism.
\end{definition}

The projection map $\pi_{\mathcal B}\colon\Sim^{\mathcal A}\to\R^d$ appearing in Birch's Theorem induced by
the choice ${\mathcal B}=\{{\bf b}_{\bf a}={\bf a}\mid{\bf a}\in{\mathcal A}\}$ of control points is called the 
\DeCo{{\sl tautological projection}} and written \DeCo{$\pi_{\mathcal A}$}.
Restricting the tautological projection to $X_{\mathcal A}$ gives the 
\DeCo{{\sl algebraic moment map}} $\DeCo{\mu}\colon X_{\mathcal A}\to\Delta$.
The components of its inverse $\mu^{-1}\colon\Delta\xrightarrow{\ \sim\ }X_{\mathcal A}$ 
provide a preferred set of blending functions for the patch.
When ${\mathcal A}\subset\Z^d$, these were studied in~\cite{CGS-GS,CGS-So03}, where they were shown to
have linear precision, and that they may be computed by 
\DeCo{{\sl iterative proportional fitting (IPF)}}, a numerical algorithm from
statistics~\cite{CGS-DR72}. 
These same arguments apply to 
irrational toric patches---the preferred
blending functions have linear precision and are computed by IPF.

Any patch has unique blending functions with linear precision~\cite[Theorem 1.11]{CGS-GS}.
While the classification of toric patches for which these preferred blending functions are
rational functions remains open in general, it has been settled for surface patches
($d=2$)~\cite{CGS-vBRS}, and this places very strong restrictions on higher-dimensional
patches.

%
%
\subsection{Iterative Proportional Fitting for Toric Patches}

In algebraic statistics, $\Sim^{\mathcal A}$ is identified with the probability simplex 
parameterizing probability distributions on data indexed by ${\mathcal A}$. 
The image $\DeCo{X_{{\mathcal A},w}}:=w.X_{\mathcal A}$ of $\R^{d}_{>}$ under the map
$\varphi_{{\mathcal A}}$~\eqref{CGS-Eq:monom_parametrize} and translation by $w$~\eqref{CGS-Eq:w-map}
is known as a \DeCo{{\sl toric model}}~\cite[\S 1.2]{CGS-PS05}.
It is more common to call this a \DeCo{{\sl log-linear model}}, as the logarithms of the 
coordinates of $\varphi_{{\mathcal A}}$ are linear functions in the logarithms of the 
coordinates of $\R^{d}_{>}$, or a \DeCo{{\sl discrete exponential family}} as 
the coordinates of  $\varphi_{{\mathcal A}}$ are exponentials in the logarithms of the 
coordinates of $\R^{d}_{>}$.

The tautological map appears in statistics as follows.
Given (observed) normalized data $q \in \Sim^{\mathcal A}$, the problem of 
\DeCo{{\sl maximum likelihood estimation}} asks for a probability distribution in the
toric model, $p\in X_{{\mathcal A},w}$, with the same sufficient statistics as $q$, 
$\mu(p) = \pi_{{\mathcal A}}(p) = \pi_{{\mathcal A}}(q)$.
By Birch's Theorem, the point $p\in X_{{\mathcal A},w}$ is unique and hence
\[
    p\ =\ \mu^{-1}(\pi_{{\mathcal A}}(q))\enspace.
\]
Thus inverting the tautological projection is necessary for maximum likelihood
estimation.  

Darroch and Ratcliff~\cite{CGS-DR72} introduced the numerical algorithm of iterative
proportional fitting, also known as generalized iterative scaling, 
for computing the inverse $\mu^{-1}$ of the tautological projection.
We now describe their algorithm.

Observe first that the toric patch $X_{{\mathcal A},w}$ does not change if we translate 
all elements of ${\mathcal A}$ by a fixed vector ${\bf b}$, (${\bf a}\mapsto {\bf a}+{\bf b}$),
so we may assume that ${\mathcal A}$ lies in the positive orthant $\R^d_>$.
Scaling the exponent vectors in ${\mathcal A}$ by a fixed positive scalar $t\in\R_>$  also does
not change $X_{{\mathcal A},w}$ as $x\mapsto x^t$ is a homeomorphism of $\R_>$ which 
extends to a homeomorphism of $\R^d_>$.
Thus we may assume that ${\mathcal A}$ lies in the standard simplex $\Delta_d$ in $\R^d$,
\[
  \Delta_d\ =\ \{ x\in\R^d_\geq \mid  {\textstyle \sum x}\leq 1\}\enspace.
\]
Lastly, we lift this to the probability simplex 
$\Sim^{d+1} \subset\R^{d+1}_{\geq}$,
\[ 
   \bSim^{d+1}\ :=\  \{ y\in \R^{d+1}_\geq\mid {\textstyle \sum y}=1\}\enspace,
\]
by 
\[
   {\mathcal A}\ni {\bf a} \ \longmapsto\ {\bf a}^+:=(1{-}{\textstyle\sum{\bf a}},\,{\bf a}) \in \Sim^{d+1}\enspace.
\]
Since for $t\in\R_>$ and $x\in\R^d_>$,
\[
   (t,tx)^{{\bf a}^+}\ =\ t^{1-\sum{\bf a}}(tx)^{\bf a}\ =\ 
    t x^{\bf a}\enspace,
\]
we see that replacing ${\mathcal A}$ by this homogeneous version also does not change 
$X_{{\mathcal A}}$.

We describe the algorithm of \DeCo{{\sl iterative proportional fitting}}, which is Theorem~1  
in~\cite{CGS-DR72}. 

\begin{proposition}\label{CGS-P:IPF}
 Suppose that ${\mathcal A}\subset\Sim^{d+1}$ has convex hull $\Delta$
 and $q\in\Sim^{\mathcal A}$.
 Set $y:=\pi_{\mathcal A}(q)\in\Delta$.
 Then the sequence of points
\[
   \{ p^{(m)}\;\mid\; m=0,1,2\dotsc\} \ \subset\ \Sim^{\mathcal A}
\]
 whose ${\bf a}$-coordinates are defined by $p_{\bf a}^{(0)}:=w_{\bf a}$ and, for 
 $m\geq 0$, 
\[
   p_{\bf a}^{(m+1)}\ :=\ p_{\bf a}^{(m)}\cdot \frac{y^{\bf a}}{(\pi_{{\mathcal A}}(p^{(m)}))^{\bf a}}\enspace,
\]
 converges to the unique point $p\in X_{{\mathcal A},w}$ such that 
 $\mu(p)=\pi_{\mathcal A}(p)=\pi_{{\mathcal A}}(q)=y$.
\end{proposition}

We remark that if ${\mathcal A}$ is not homogenized then to compute 
$\mu^{-1}(y)$ for $y\in \Delta$, we first put
${\mathcal A}$ into homogeneous form using an affine map $\psi$, and then 
use iterative proportional fitting to compute
$\pi_{{\mathcal A}^{+}}^{-1}(\psi(y))=\pi_{{\mathcal A}}^{-1}(y)$.  
We also call this modification of the algorithm of Proposition~\ref{CGS-P:IPF}
iterative proportional fitting.
Thus iterative proportional fitting computes the inverse image of the tautological
projection.

%
\section{Injectivity of Patches}\label{CGS-S:inj}
%
Birch's Theorem (Theorem~\ref{CGS-Th:Birch}) states that for one particular choice of 
control points, namely $\{{\bf b}_{\bf a}:={\bf a}\mid{\bf a}\in{\mathcal A}\}$ and all weights 
$w\in\R^{\mathcal A}_>$,
the mapping $F$~\eqref{CGS-Eq:P_Patch} of a toric patch of shape $({\mathcal A},w)$ is a homeomorphism
onto its image.  
From this, we can infer that for {\sl most} choices of control points and weights, this
mapping is at least an immersion.
To study dynamical systems arising from chemical reaction networks, Craciun and
Feinberg~\cite{CGS-CF05} prove an injectivity theorem for certain maps, which we adapt to 
generalize Birch's Theorem.
This will give conditions on control points  ${\mathcal B}\subset\R^d$ which guarantee that for
{\sl any} choice $w$ of weights, the resulting mapping $F$~\eqref{CGS-Eq:P_Patch} of a toric
patch of shape $({\mathcal A},w)$ is a homeomorphism onto its image.
This result has several consequences concerning the injectivity of toric patches.

Let us first give the Craciun-Feinberg Theorem.
Let $Y=\{{\bf y}_1,\dotsc,{\bf y}_m\}\subset\R^n$ be a finite set of points which affinely
spans $\R^n$. 
For $k\in\R^m_>$ and $Z:=\{{\bf z}_1,\dotsc, {\bf z}_m\}\subset\R^n$,
consider the map $\DeCo{G_k}\colon\R^n_>\to\R^n$ defined by
 \begin{equation}\label{CGS-Eq:exp_sum}
    G_k(x)\ :=\ \sum_{i=1}^m k_i\, x^{{\bf y}_i} {\bf z}_i\enspace.
 \end{equation}
%

\begin{theorem}[Craciun-Feinberg]\label{CGS-Th:CF}
  The map $G_k$ is injective for every $k\in\R^m_>$ if and only 
  if the determinant of the Jacobian matrix,
\[
   \mbox{\rm Jac}(G_k)\ =\ 
    \left(\frac{\partial (G_k)_i}{\partial x_j}\right)_{i,j=1}^n\enspace ,
\]
  does not vanish for any $x\in\R^n_>$ and any $k\in\R^m_>$.
\end{theorem}

We give a proof in Appendix~\ref{CGS-A:A}.

The condition of Theorem~\ref{CGS-Th:CF} that the Jacobian $\mbox{\rm Jac}(G_k)$ does not
vanish for any $x\in\R^n_>$ is reminiscent of the Jacobian conjecture~\cite{CGS-Ke}, 
which is that the Jacobian of a polynomial map $G\colon\C^n\to\C^n$ does not vanish if and
only if the map $G$ is an isomorphism.
Since we are restricted to $x\in\R^n_>$, it is closer to the real Jacobian
conjecture, which is however false~\cite{CGS-Pi}, and therefore not necessarily relevant.

The condition of Theorem~\ref{CGS-Th:CF} is conveniently restated in terms of $Y$ and $Z$.
For a list $I=\{i_1,\dotsc,i_n\}\subset\{1,\dotsc,m\}$, which we write as
$I\in\DeCo{\binom{[m]}{n}}$, let $Y_I$ be the determinant of the matrix whose columns are
the vectors ${\bf y}_{i_1},\dotsc,{\bf y}_{i_n}$, and define $Z_I$ similarly.
In Appendix~\ref{CGS-A:A}, we deduce the following corollary.

\begin{corollary}\label{CGS-C:matroid}
 The map $G_k$~$\eqref{CGS-Eq:exp_sum}$ is injective for all $k\in\R^m_>$ 
 if and only if $(Y_I\cdot Z_I)\cdot(Y_J\cdot Z_J)\geq 0$ for every
 $I,J\in\binom{[m]}{n}$ and at least one product $Y_I\cdot Z_I$ is non-zero.
\end{corollary}

This leads to a generalization of Birch's Theorem.
An ordered list $p_0,\dotsc,p_d$ of affinely independent points in $\R^d$
determines an orientation of $\R^d$---simply consider the basis
\[
    p_1{-}p_0\,,\ p_2{-}p_0\,,\ \dotsc\,,\ p_d{-}p_0\enspace.
\]
Let ${\mathcal A}$ and ${\mathcal B}=\{{\bf b}_{\bf a}\mid{\bf a}\in{\mathcal A}\}$ be finite
sets of points in $\R^d$. 
Suppose that $\{{\bf a}_0,\dotsc,{\bf a}_d\}$ is an affinely independent subset
of ${\mathcal A}$.
If the corresponding subset $\{{\bf b}_{{\bf a}_0},\dotsc,{\bf b}_{{\bf a}_d}\}$ of ${\mathcal B}$ is also
affinely independent, then each  subset determines an orientation, and the
two orientations are either the same or they are opposite.
We say that ${\mathcal A}$ and ${\mathcal B}$ are  \DeCo{{\sl compatible}} if either every such pair of
orientations is the same, or if every such pair of orientations is opposite.
We further need that 
there is at least one affinely independent subset of
${\mathcal A}$ such that the corresponding subset of ${\mathcal B}$ is also affinely
independent. 
Observe that compatibility is preserved by invertible affine transformations acting
separately on ${\mathcal A}$ and ${\mathcal B}$.

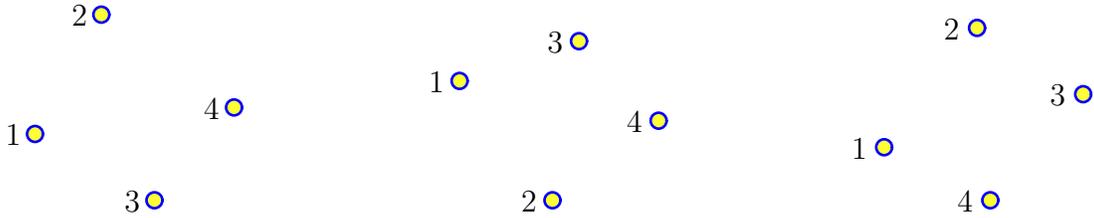
\begin{figure}[htb]
 \[
  \begin{picture}(410,85)
   \put(0,0){\includegraphics{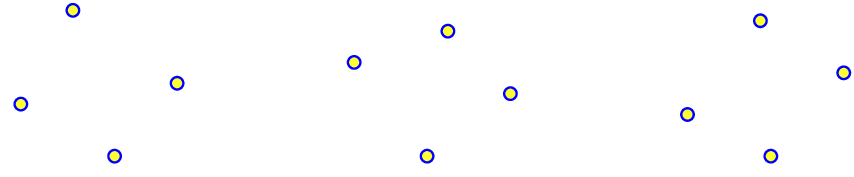}}
   \put( -1,31){$1$}
   \put( 24,76){$2$}
   \put( 44, 6){$3$}
   \put( 74,41){$4$}

   \put(159,51){$1$}
   \put(204,66){$3$}
   \put(194, 6){$2$}
   \put(234,36){$4$}

   \put(319,26){$1$}
   \put(354,71){$2$}
   \put(359, 6){$4$}
   \put(394,46){$3$}
  \end{picture}
 \]
 \caption{Compatible and incompatible sets of points.}\label{CGS-F:Compatible}
\end{figure}

In Fig.~\ref{CGS-F:Compatible} shows three sets of labeled points.
The first and second sets are compatible, but neither is compatible with the third.

We give our generalization of Birch's Theorem.
Suppose that $\Delta\subset\R^d$ is the convex hull of ${\mathcal A}$ and
$\{\beta_{\bf a}\colon\Delta\to\R_\geq\mid {\bf a}\in{\mathcal A}\}$ are toric B\'ezier functions for
${\mathcal A}$.
For any $w\in\R^{\mathcal A}_>$, let $F_w\colon \Delta\to \R^d$ be the toric patch of shape
$({\mathcal A},w)$ given by the control points ${\mathcal B}\subset\R^d$:
 \begin{equation}\label{CGS-Eq:T_Patch}
   F_w(x)\ :=\ 
   \frac{\sum_{{\bf a}\in{\mathcal A}} w_{\bf a}\beta_{\bf a}(x)\,{\bf b}_{\bf a} }
        {\sum_{{\bf a}\in{\mathcal A}} w_{\bf a}\beta_{\bf a}(x)}\enspace .
 \end{equation}
%

\begin{theorem}\label{CGS-Th:new_birch}
 The map $F_w$  is injective  for all $w\in\R^{\mathcal A}_>$ if and only if 
 ${\mathcal A}$ and ${\mathcal B}$ are compatible.
\end{theorem}

As any set ${\mathcal A}$ is compatible with itself, this implies Birch's
Theorem (Theorem~\ref{CGS-Th:Birch}).

\begin{example}\label{CGS-Ex:cubic_move}
 Let $\btri$ be the convex hull of $\{(0,0),(1,0),(0,1)\}$.
 Set ${\mathcal A}:=3\,\tri\cap\Z^2$ and let $w\in\R^{\mathcal A}_>$
 be the weights of a cubic B\'ezier patch (Example~\ref{CGS-Ex:BezTri} with $m=3$).
 We consider choices ${\mathcal B}\subset\R^2$ of control points that are compatible with
 ${\mathcal A}$. 
 For convenience, we will require that ${\bf b}_{\bf a}={\bf a}$ when ${\bf a}$ is a vertex and 
 that if ${\bf a}$ lies on an edge of $3\,\tri$, then so does ${\bf b}_{\bf a}$.
 For these edge control points, compatibility imposes the restriction that they appear
 along the edge in the same order as the corresponding exponents from ${\mathcal A}$.
 The placement of the center control point is however constrained.
 We show two compatible choices of ${\mathcal B}$ in Fig.~\ref{CGS-F:Compatible_Control}.
 On the left is the situation of Birch's Theorem, in which ${\bf b}_{\bf a}={\bf a}$, and on the right
 we have moved the edge control points.
 The region in which we are free to move the center point is shaded in each picture.
 \begin{figure}[htb]
\[
  \includegraphics{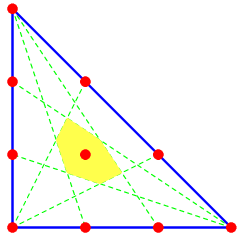}
   \qquad\qquad
   \includegraphics{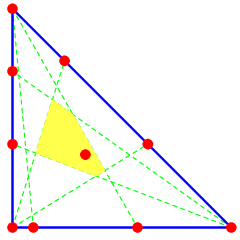}
\]
 \caption{Compatible control points for the B\'ezier triangle.}\label{CGS-F:Compatible_Control}
\end{figure}
\end{example}

\begin{proof}[Theorem~$\ref{CGS-Th:new_birch}$]
 Let $(t,x)$ be coordinates for $\R^{d+1}$ and consider the map 
 $G_w\colon\R^{d+1}_>\to \R^{d+1}$ defined by
\[
   G_w(t,x)\ =\ \sum_{{\bf a}\in{\mathcal A}} t x^{{\bf a}}w_{\bf a}(1,{\bf b}_{\bf a})\enspace.
\]
 We claim that $F_w$ is injective if and only if $G_w$ is injective.

 Since $F_w$ is the composition~\eqref{CGS-Eq:composition}
 $\Delta\xrightarrow{\,\beta\,}X_{\mathcal A}\xrightarrow{\,w.\,}w.X_{\mathcal A}\xrightarrow{\,\pi_{\mathcal B}\,}\R^d$,
 with the first  map an isomorphism, $F_w$ is injective if and only if the 
 composition of the last two maps is injective.
 Since $X_{\mathcal A}$ is compact, this will be injective if and only if its restriction to
 the interior $X^\circ_{\mathcal A}$ of $X_{\mathcal A}$ is injective.
 Precomposing with the monomial parametrization~\eqref{CGS-Eq:monom_parametrize} of
 $X^\circ_{\mathcal A}$, we see that $F_w$ is injective if and only if the map $H_w\colon
 \R^d_>\to \R^d$ defined by 
\[
   H_w(x)\ =\ 
   \frac{\sum_{{\bf a}\in{\mathcal A}} x^{\bf a} w_{\bf a} {\bf b}_{\bf a}}{\sum_{{\bf a}\in{\mathcal A}} x^{\bf a} w_{\bf a}}
   \ \colon\ \R^d_>\ \xrightarrow{\,\varphi_{\mathcal A}\,}\ 
    X_{\mathcal A}\ \xrightarrow{\,w.\,}\ w.X_{\mathcal A}\ 
    \xrightarrow{\,\pi_{\mathcal B}\,}\ \R^d
\]
 is injective.

 Since $G_w(t,x)=t\cdot G_w(1,x)$, these values lie on a ray through the origin, and we
 invite the reader to check that this ray meets the hyperplane with first coordinate 1 at
 the point $(1,H_w(x))$.  
 Thus $G_w$ is injective if and only if $H_w$ is injective, which is equivalent to the
 map $F_w$ being injective.

 We deduce the theorem by showing that $G_w$ is injective.
 This follows from Corollary~\ref{CGS-C:matroid} 
 as the condition that ${\mathcal A}$ and ${\mathcal B}$ are compatible is
 equivalent to $Y_I\cdot Z_I\geq 0$ for all $I\in\binom{[m]}{d+1}$, where
 $Y=\{(1,{\bf a})\mid{\bf a}\in{\mathcal A}\}$, 
 $Z=\{(1,{\bf b}_{\bf a})\mid {\bf a}\in{\mathcal A}\}$ and we have $m=\#{\mathcal A}$.
\end{proof}

We now describe two applications of Theorem~\ref{CGS-Th:new_birch} to modeling.

%
%
\subsection{Solid Modeling with (Toric) B\'ezier Patches}
In solid modeling, we represent a 3-dimensional solid by covering it with 3-dimensional
patches, for example using B\'ezier toric patches as finite elements.
Besides the obvious $C^0$ or higher continuity along the boundary as required, such
B\'ezier finite elements should at least provide a one-to-one parametrization of their
image, i.e. they should be injective.
By Theorem~\ref{CGS-Th:new_birch}, we may guarantee injectivity by requiring that the control
points ${\mathcal B}$ be compatible with the exponents ${\mathcal A}$.
Moreover, if these sets are incompatible, then there is some choice of weights for which
the patch is not injective.

%
%
\subsection{Injectivity of B\'ezier Curves and Surfaces}
Typically, the exponents ${\mathcal A}$ and the control points do not lie in the same
space; surfaces (${\mathcal A}\subset\R^2$) are modeled in 3-space (${\mathcal B}\subset \R^3$),
or curves (${\mathcal A}\subset\R$) in 2-- or 3--space 
(${\mathcal B}\subset \R^2\mbox{\ or\ }\R^3$).
Nevertheless, Theorem~\ref{CGS-Th:new_birch} gives conditions that imply injectivity of
patches.

Let $p\in\R^{n+1}$ be a point disjoint from a hyperplane, $H$.
The projection $\R^{n+1}-\to H$ with center $p$ is the map which associates a
point $x\in\R^{n+1}$ to the intersection of the line $\overline{px}$ with $H$.
We use a broken arrow as the projection is not defined on the plane through $p$ parallel to $H$.
Identifying $H$ with $\R^n$ gives a projection map $\R^{n+1}-\to\R^n$.
A coordinate projection $(x_1,\dotsc,x_n,x_{n+1})\mapsto(x_1,\dotsc,x_n)$ is a projection
with center at infinity.
More generally, a \DeCo{{\sl projection}} $\R^n-\to\R^d$ is a sequence of such projections from points.

\begin{theorem}\label{CGS-Th:inj}
 Let ${\mathcal A}\subset\R^d$, $w\in\R^{\mathcal A}_>$, and 
 ${\mathcal B}\subset\R^n$ be the exponents, weights, 
 and control points of a toric patch and let $\Delta$ be the convex hull of ${\mathcal A}$.
 If there is a projection $\pi\colon\R^n-\to\R^d$ such that ${\mathcal A}$ is compatible with the
 image $\pi({\mathcal B})$ of ${\mathcal B}$, then the mapping $F\colon\Delta\to\R^n$ given by the toric
 blending functions associated to ${\mathcal A}$, the weights $w$, and control points ${\mathcal B}$ is
 injective. 
\end{theorem}

\begin{proof}
 By Theorem~\ref{CGS-Th:new_birch}, the composition $\pi\circ F$ is injective, from which it
 follows that $F$ must have been injective.
\end{proof}
 
\begin{example}
  For the curve on the left in Fig.~\ref{CGS-F:one} (which is reproduced below), the
  vertical projection 
  $\R^2\to\R^1$ maps the control points $\{{\bf b}_0,\dotsc,{\bf b}_5\}$
  to points on the line in the same order as the exponents ${\mathcal A}=\{0,1,2,3,4,5\}$,
  which implies that the curve has no self-intersections.
 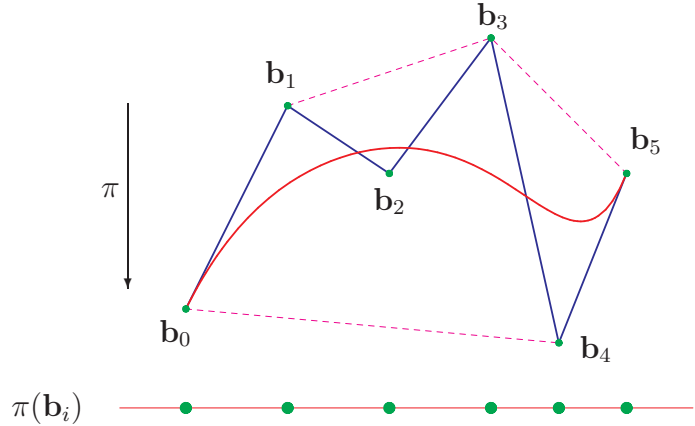
\begin{figure}[htb]
 \[
  \begin{picture}(270,171)(-64,-30)
     \put(-26,-19){\includegraphics{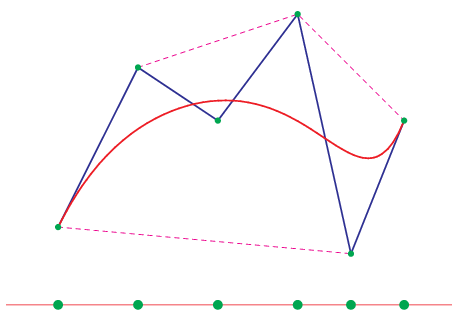}}
     \put(-30,65){$\pi$}\put(-20,100){\vector(0,-1){70}}
     \put(170,83){$\bb_5$}
     \put(151,4){$\bb_4$}
     \put(112,130){$\bb_3$}
     \put(73,60){$\bb_2$}
     \put(32,108){$\bb_1$}
     \put(-8,10){$\bb_0$}
     \put(-64,-18){$\pi(\bb_i)$}  
   \end{picture}
 \]
 \caption{A compatible projection.}\label{CGS-F:Projection_compatible}
\end{figure}
\end{example}
 
%
\section{Control Polytopes and Toric Degenerations}\label{CGS-S:degeneration}
%
The convex hull property asserts that the image, $F(\Delta)$, of a toric B\'ezier patch
of shape $({\mathcal A},w)$ given by control points ${\mathcal B}=\{{\bf b}_{\bf a} \mid {\bf a}\in{\mathcal A}\}\subset\R^n$
and weights $w\in\R^{\mathcal A}_>$ lies in the convex hull of the control points.
When $F(\Delta)$ is a curve, the control points may be joined sequentially to form the
\DeCo{{\sl control polygon}}, which is a piecewise linear representation of the curve.
When $F(\Delta)$ is however a surface patch, there are many ways to interpolate
the control points by triangles or other polygons to obtain a piecewise linear surface,
called a \DeCo{{\sl control polytope}}, that represents the patch.
The shape of this control polytope affects the shape of the patch.
For example, when the control points have the form $({\bf a},\lambda({\bf a}))$ for $\lambda$ a
convex function, then the patch is convex~\cite{CGS-DM88,CGS-Da91}.
Also, Leroy~\cite{CGS-Leroy} uses a particular control polytope for the graph of a
function to obtain certificates of positivity for polynomials.

Among all control polytopes for a given set of control points, we identify the class of
regular control polytopes, which come from regular triangulations of the exponents
${\mathcal A}$.
These regular control polytopes are related to the shape of the patch in the following
precise manner: 
There is a choice of weights so that a toric B\'ezier patch is arbitrarily close to a
given control polytope if and only if that polytope is regular.

%
\subsection{B\'ezier Curves}
%
It is instructive to begin with B\'ezier curves.
A B\'ezier curve of degree $m$ in $\R^n$ with weights $w$ is the 
composition~\eqref{CGS-Eq:composition},  
\[
   [0,1]\ \xrightarrow{\;\beta\;}\ \Sim^{m+1}\ \xrightarrow{\;w.\;}\ 
    \Sim^{m+1}\ \xrightarrow{\;\pi_{\mathcal B}\;}\ \R^n\enspace,
\]
where $\beta=(\beta_0,\dotsc,\beta_m)$ with 
$\beta_i(x)=\binom{m}{i}x^i(1-x)^{m-i}$ for $x\in[0,1]$.
Then the map $\beta$ is given by $z_i=\binom{m}{i}x^i(1-x)^{m-i}$,
for $i=0,\dotsc,m$.
Here, $(z_0,\dotsc,z_m)\in\R^{d+1}_\geq$ with $z_0+\dotsb+z_m=1$
are the coordinates for $\Sim^{m+1}$.
The image $\beta[0,1]\subset\Sim^{m+1}$ is defined by the binomials
 \begin{equation}\label{CGS-Eq:Toric_binom}
   \tbinom{m}{i} \tbinom{m}{j} z_a z_b 
    \ -\ 
   \tbinom{m}{a}\tbinom{m}{b} z_i z_j
   \ =\ 0\,,\qquad\mbox{for}\quad a+b=i+j\enspace.
 \end{equation}
To see this, suppose that $(z_0,\dotsc,z_m)\in\R^{m+1}_\geq$ satisfies~\eqref{CGS-Eq:Toric_binom}.
 Setting $x:=z_1/(mz_0+z_1)$, then we may solve these equations to obtain 
 $z_i=\binom{m}{i}x^i(1-x)^{m-i}=\beta_i(x)$.

If $w=(w_0,\dotsc,w_m)\in\R^{m+1}_>$ are weights, then 
$w.\beta[0,1]$ is defined in $\Sim^{m+1}$ by
 \begin{equation}\label{CGS-Eq:complicated}
   w_iw_j \tbinom{m}{i}\tbinom{m}{j} z_a z_b\ -\ 
   w_aw_b \tbinom{m}{a}\tbinom{m}{b} z_i z_j
   \ =\ 0\,,\qquad\mbox{for}\quad a+b=i+j\enspace.
 \end{equation}
Suppose that we choose weights $\DeCo{w_i}:=t^{i(m-i)}$.
Dividing by $t^{i(m-i)+j(m-j)}$,~\eqref{CGS-Eq:complicated} becomes
 \begin{equation}\label{CGS-Eq:binoms}
   \tbinom{m}{i}\tbinom{m}{j} z_a z_b\ -\ 
   t^{i^2+j^2-a^2-b^2} \tbinom{m}{a}\tbinom{m}{b} z_i z_j
   \ =\ 0\,,\qquad\mbox{for}\quad a+b=i+j\enspace.
 \end{equation}
Since $i+j=a+b$, we may assume that $a<i\leq j<b$.
Setting $\DeCo{c}:=i-a=b-j\geq 1$, we see that
\[
  i^2+j^2-a^2-b^2\ =\  c(i+a)-c(j+b)\ =\ c(i-j+a-b)\leq -2c\ <\ 0\enspace.
\]
If we consider the limit of these binomials~\eqref{CGS-Eq:binoms} as $t\to\infty$, 
we obtain
\[
    z_a\cdot z_b\ =\ 0\qquad\mbox{if}\quad |a-b|>1\enspace.
\]
These define the polygonal path in $\Sim^{m+1}$ whose $i$th segment is the edge 
\[
  ( \underbrace{0,\dotsc,0}_{i-1}\,,\,x\,,\, 1-x\,,\,
    \underbrace{0,\dotsc,0}_{m-i})\qquad\mbox{\rm for}\quad x\in[0,1]\enspace,
\]
and whose projection to $\R^n$ is the \DeCo{{\sl control polygon}} of the B\'ezier curve, 
which is the collection of line segments $\overline{{\bf b}_0,{\bf b}_1}$,
$\overline{{\bf b}_1,{\bf b}_2}$, $\dotsc$, $\overline{{\bf b}_{m-1},{\bf b}_m}$.

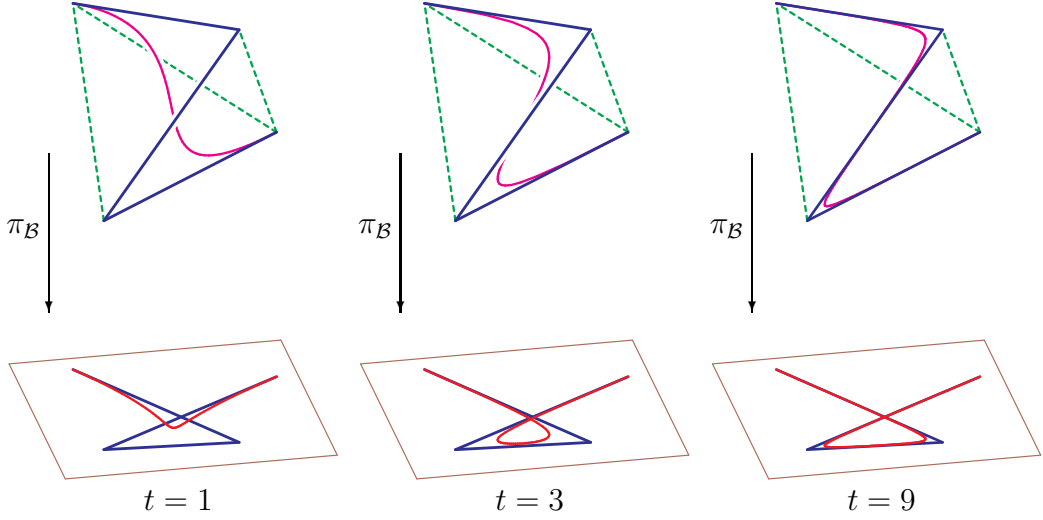
\begin{figure}[htb]
\[
  \begin{picture}(135,192)(-2,0)
    \put(-1,10){\includegraphics[height=185pt]{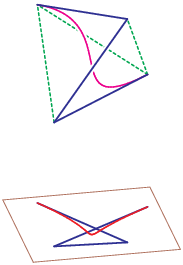}}
    \put(52,0){$t=1$}
    \put(0,105){$\pi_{\mathcal B}$} \put(16,135){\vector(0,-1){60}}
  \end{picture}
  \begin{picture}(133,187)
    \put(-1,10){\includegraphics[height=185pt]{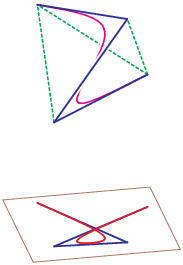}}
    \put(52,0){$t=3$}
    \put(0,105){$\pi_{\mathcal B}$} \put(16,135){\vector(0,-1){60}}
   \end{picture}
  \begin{picture}(133,187)
    \put(-1,10){\includegraphics[height=185pt]{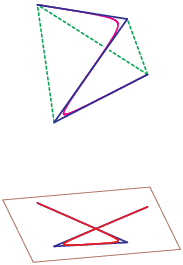}}
    \put(52,0){$t=9$}
    \put(0,105){$\pi_{\mathcal B}$} \put(16,135){\vector(0,-1){60}}
  \end{picture}
\]
\caption{Toric degenerations of a B\'ezier cubic.}
\label{CGS-F:cubic_bezier}
\end{figure}

We illustrate this when $m=3$ in Fig.~\ref{CGS-F:cubic_bezier},
which shows three different B\'ezier curves having the same  control points, but different
weights $w_i=t^{i(3-i)}$ for $t=1,3,9$.
In algebraic geometry, altering the weights in this manner is called a 
\DeCo{{\sl toric degeneration}}.
The B\'ezier cubics are displayed together with the cubics $w.\beta[0,1]$ lying in
the 3-simplex, $\includegraphics{figures/Simplex.eps}^4$, which is drawn in $\R^3$.
In these pictures, the projection $\pi_{\mathcal B}$ is simply the vertical projection forgetting
the third coordinate.
The progression indicated in Fig.~\ref{CGS-F:cubic_bezier}, where the B\'ezier curve
approaches the control polygon as the parameter $t$ increases, is a general phenomenon.
Let $\|\cdot\|$ be the usual Euclidean distance in $\R^n$.

\begin{theorem}\label{CGS-Th:Bez_curve_approx}
 Suppose that $F_t\colon[0,1]\to\R^n$ is a B\'ezier curve of degree $m$ with 
 control points ${\mathcal B}$ and weights $w_i=t^{i(m-i)}$.
 Set $\DeCo{\kappa}:=\max\{\|{\bf b}_{\bf a}\|: {\bf b}_{\bf a}\in{\mathcal B}\}$.
 For any $\epsilon>0$, if we have $t>\kappa m/\epsilon$, then the distance between the control
 polygon and any point of the B\'ezier curve $F_t[0,1]$ is less than $\epsilon$. 
\end{theorem}

\begin{proof}
 Let $z\in w.\beta[0,1]\subset\Sim^{m+1}$, where the weights are $w_i=t^{i(m-i)}$ with
 $t>\kappa m/\epsilon$. 
 Suppose that $b-a>1$ are integers in $[0,m]$.
 Then there exist integers $i\leq j$ with $a<i\leq j<b$ and $a+b=i+j$.
 By~\eqref{CGS-Eq:binoms}, 
 we have
\[
   z_az_b\ =\ t^{i^2+j^2-a^2-b^2}
   \frac{\tbinom{m}{a}\tbinom{m}{b}}{\tbinom{m}{i}\tbinom{m}{j}} z_iz_j\enspace .
\]
Since the binomial coefficients are log-concave,
we have
\[
   \binom{m}{a}\binom{m}{b}\ <\ 
   \binom{m}{i}\binom{m}{j}\enspace .
\]
Using $i^2+j^2-a^2-b^2<-2$ and $z_i+z_j\leq 1$, we see that 
\[
   z_az_b\ <\ \frac{1}{4t^2}\enspace .
\]

In particular, if $|b-a|>1$, then at most one of $z_a$ or $z_b$ exceeds $1/2t$.
We conclude that at most two, necessarily consecutive, coordinates of $z$ may exceed
$1/2t$.  
Suppose that $i$ is an index such that $z_j<1/2t$ if $j\neq i{-}1,i$
and let $\DeCo{x}:= z_{i-1}{\bf b}_{i-1} + (1-z_{i-1}){\bf b}_i$, a point along the $i$th
segment of the control polygon. 
Since 
\[
  1\ \geq\ z_{i-1}+z_i\ =\ 1-\sum_{j\neq i-1,i} z_j\ >\ 
  1-\frac{m-1}{2t}\enspace ,
\]
we have $0\leq1-z_{i-1}-z_i<\frac{m-1}{2t}<\frac{m}{2t}$, and we see that 
 \begin{eqnarray*}
    \left\| \pi_{\mathcal B}(z) - x\right\| &=&
    \left\| \sum z_j{\bf b}_j\ -\ \left( z_{i-1}{\bf b}_{i-1} + (1-z_{i-1}){\bf b}_i\right)\right\|\\
    &\leq&
    \sum_{j\neq i-1,i} z_j \|{\bf b}_j\|\ +\ 
    |z_i-(1-z_{i-1})|\|{\bf b}_i\|\\
    &<&  \kappa\,\frac{m-1}{2t} \ +\ \kappa\,\frac{m}{2t}\ <\ \frac{\kappa m}{t}\ =\ \epsilon\,.
 \end{eqnarray*}
This proves the theorem as $z\in w.\beta[0,1]$ is an arbitrary point of the curve $F[0,1]$.
\end{proof}

%
\subsection{Regular Triangulations and Control Polytopes}\label{CGS-S:regular}
%
In dimensions $2$ and higher, the analog of Theorem~\ref{CGS-Th:Bez_curve_approx} requires the
notion of a regular triangulation from geometric combinatorics.
Let ${\mathcal A}\subset\R^d$ be a finite set of points and consider a lifting function
$\lambda\colon{\mathcal A}\to\R$. 
Let \DeCo{$P_\lambda$} be the convex hull of the lifted points 
\[
   \lambda({\mathcal A})\ :=\ 
   \{ ({\bf a},\,\lambda({\bf a})) \mid {\bf a}\in{\mathcal A}\}\ 
    \subset\ \R^{d+1}\enspace.
\]
We assume that $P_\lambda$ is full-dimensional in that $\R^{d+1}$ is its affine span.

The \DeCo{{\sl upper facets}} of $P_\lambda$ are those facets whose outward pointing normal
vector has positive last coordinate.
Any face of an upper facet is an \DeCo{{\sl upper face}}.
We illustrate this below when $d=1$, where the displayed arrows are outward pointing
normal vectors to upper facets. 
%
 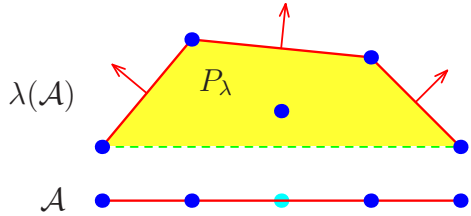
\begin{figure}[htb]
 \[
  \begin{picture}(180,75)(-32,0)
   \put(0,0){\includegraphics[height=80pt]{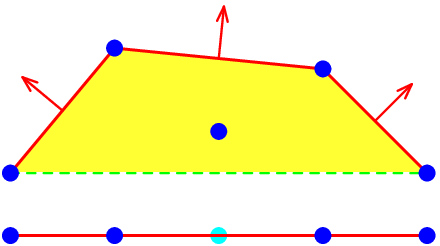}}
   \put(-32,40){$\lambda({\mathcal A})$}
   \put(40,45){$P_\lambda$}
   \put(-20,0){${\mathcal A}$}
  \end{picture}
 \]
 \caption{Upper faces and a regular subdivision.}\label{CGS-F:Upper_face}
\end{figure}

Projecting these upper facets to $\R^d$ yields a \DeCo{{\sl regular polyhedral subdivision}}
of the convex hull of ${\mathcal A}$, which is the image of $P_\lambda$.
For our purposes, we will need to assume that the lifting function is generic in that all
upper facets of $P_\lambda$ are simplices.
In this case, we obtain a \DeCo{{\sl regular triangulation}} of ${\mathcal A}$.
This consists of a collection 
\[
   \{ {\mathcal A}_i\mid i=1,\dotsc,m\}
\]
of subsets of ${\mathcal A}$, where each subset ${\mathcal A}_i$ consists of $d{+}1$ elements and
spans a $d$-dimensional simplex. 
We regard all subsets of the facets ${\mathcal A}_i$ as faces of the triangulation.
These simplices form a subdivision in that they cover the convex
hull of ${\mathcal A}$ and any two with a non-empty intersection meet along a common face.

The subdivision induced by the lifting function in Fig.~\ref{CGS-F:Upper_face} consists of
three intervals resulting from removal of the middle point of ${\mathcal A}$, which does not
participate in the subdivision as it is not lifted high enough.

\begin{figure}[htb]
\[
  \includegraphics[height=90pt]{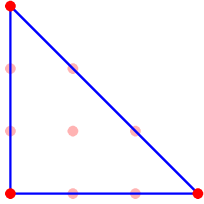}\qquad
  \includegraphics[height=90pt]{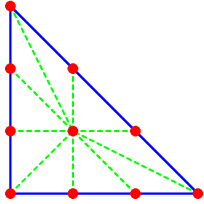}\qquad
  \includegraphics[height=90pt]{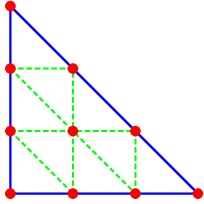}\qquad
  \includegraphics[height=90pt]{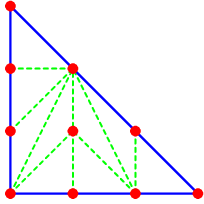}
 \]
 \caption{Some regular triangulations.}\label{CGS-F:regular_triangulations}
\end{figure}

A set may have many regular triangulations, and not every point 
needs to participate in a given triangulation.
Figure~\ref{CGS-F:regular_triangulations} shows the edges in four regular triangulations of 
$3\,\tri\cap\Z^2$.
%

\begin{figure}[htb]
\[
  \includegraphics[height=120pt]{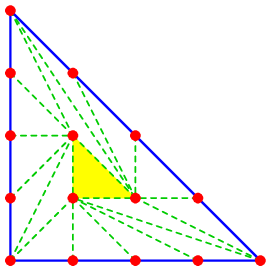}
\]
 \caption{An irregular triangulation.}\label{CGS-F:irrregular_triangulation}
\end{figure}

Not every triangulation is regular.
We may assume that a lifting function $\lambda$ for the triangulation 
of $4\,\tri\cap\Z^2$ in
Fig.~\ref{CGS-F:irrregular_triangulation} takes a constant value at the three interior points.  
The clockwise neighbor of any vertex of the big triangle must be lifted
lower than that vertex. 
(Consider the figure they form with the parallel edge of the
interior triangle.)
Since the edge of the big triangle is lifted to a convex path,
this is impossible, except in some M.C. Escher woodcuts.

\begin{definition} 
 Let ${\mathcal B} = \{{\bf b}_{\bf a} \mid {\bf a} \in {\mathcal A}\} \subset \R^n$ be 
 a collection of control points indexed by a finite set of exponents 
 ${\mathcal A} \subset \R^d$ with $d\leq n$. 
 Given a regular triangulation ${\mathcal T} = \{ {\mathcal A}_i \mid i = 1,\dots,m\}$ of 
 ${\mathcal A}$ we define the control polytope as follows. 
For each $d$-simplex ${\mathcal A}_i$ in ${\mathcal T}$,
the corresponding points of ${\mathcal B}$ span a (possibly degenerate) simplex 
\[
   \mbox{conv}\{{\bf b}_{\bf a}\mid {\bf a}\in{\mathcal A}_i\}\enspace.
\]
The union of these simplices in $\R^n$ forms the
\DeCo{{\sl regular control polytope} ${\mathcal B}({\mathcal T})$} that is
induced by the regular triangulation ${\mathcal T}$ of ${\mathcal A}$.
This is a simplicial complex in $\R^n$ with vertices in ${\mathcal B}$ that has the same
combinatorial type as the triangulation ${\mathcal T}$ of ${\mathcal A}$.
\end{definition} 

If the coordinate points $({\bf e}_{\bf a}\mid{\bf a}\in{\mathcal A})$ of $\R^{\mathcal A}$ are
our control points 
(these are the vertices of $\Sim^{\mathcal A}$), then the regular control polytope is just the 
geometric realization \DeCo{$|{\mathcal T}|$} of the simplicial complex ${\mathcal T}$, which is a 
subcomplex of the simplex $\Sim^{\mathcal A}$. 
In general, ${\mathcal B}({\mathcal T})$ is the image of this geometric realization
$|{\mathcal T}|\subset\Sim^{\mathcal A}$ under the projection $\pi_{\mathcal B}$.

\begin{example}\label{CGS-Ex:control_polytope}
  Let ${\mathcal A}:=3\,\tri\cap\Z^2$, the exponents for a cubic B\'ezier triangle.
  Figure~\ref{CGS-F:Control_polytopes} shows the three control polytopes corresponding to the last
  three regular triangulations of Figure~\ref{CGS-F:regular_triangulations}, all with the same control
  points.  
 \begin{figure}[htb]
\[
   \includegraphics[height=95pt]{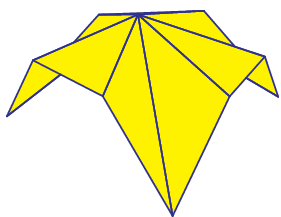}\quad
   \includegraphics[height=95pt]{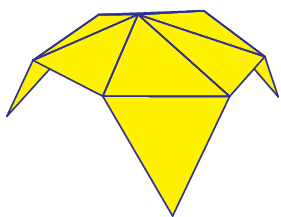}\quad
   \includegraphics[height=95pt]{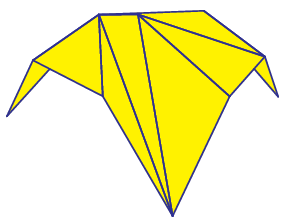}
\]
 \caption{Three control polytopes.}\label{CGS-F:Control_polytopes}
\end{figure}
\end{example}

The reason that we introduce regular control polytopes is that they may be
approximated by toric B\'ezier patches.

\begin{theorem}\label{CGS-Th:Bez_approx}
  Let ${\mathcal A}\subset\R^d$,  $w\in\R^{\mathcal A}_>$, and 
  ${\mathcal B}\subset\R^n$  be exponents,
  weights, and control points for a toric B\'ezier patch.
  Suppose that ${\mathcal T}$ is a regular triangulation of ${\mathcal A}$ induced by a lifting
  function $\lambda\colon{\mathcal A}\to\R$. 
  For each $t>0$, let $F_t\colon \Delta\to\R^n$ be the toric B\'ezier patch of shape
  ${\mathcal A}$ with control points ${\mathcal B}$ and weights $t^{\lambda({\bf a})}w_{\bf a}$.
  Then, for any $\epsilon>0$ there exists a $t_0$ such that if $t>t_0$, the image 
  $F_t(\Delta)$ lies within $\epsilon$ of the control polytope ${\mathcal B}({\mathcal T})$.
\end{theorem}

We prove Theorem~\ref{CGS-Th:Bez_approx} in Appendix~\ref{CGS-A:B}.
Figure~\ref{CGS-F:Control_polytope_degeneration} illustrates Theorem~\ref{CGS-Th:Bez_approx} for a
cubic B\'ezier triangle with the control points of Example~\ref{CGS-Ex:control_polytope}.
The patch on the left is the cubic B\'ezier triangle with the weights of
Example~\ref{CGS-Ex:BezTri}. 
The second and third patches are its deformations corresponding to the lifting function
inducing the leftmost control polytope of Fig.~\ref{CGS-F:Control_polytopes}.
The values of $t$ are $1$, $5$, and $100$, as we move from left to right.
 \begin{figure}[htb]
\[
   \includegraphics[height=95pt]{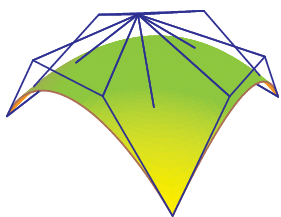}\quad
   \includegraphics[height=95pt]{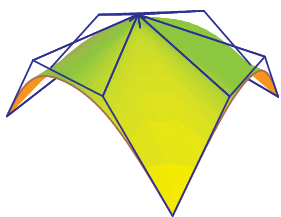}\quad
   \includegraphics[height=95pt]{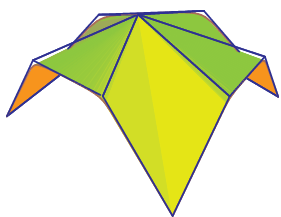}
\]
 \caption{Degeneration to the control polytope.}\label{CGS-F:Control_polytope_degeneration}
\end{figure}

An absolutely unpractical consequence of Theorem~\ref{CGS-Th:Bez_approx} is a 
universality result:
Any surface which admits a triangulation that forms a regular control polytope may
be approximated by a single B\'ezier patch.

As with Theorem~\ref{CGS-Th:Bez_curve_approx}, the main idea behind the proof of
Theorem~\ref{CGS-Th:Bez_approx} (which is given in Appendix~\ref{CGS-A:B}) is that 
for $t$ large enough, the translated patch $t.X_{{\mathcal A},w}\subset\Sim^{\mathcal A}$ can be made
arbitrarily close to the geometric realization $|{\mathcal T}|\subset\Sim^{\mathcal A}$ of the regular
triangulation ${\mathcal T}$. 
The result follows by projecting this into $\R^n$ using $\pi_{\mathcal B}$.

In Appendix~\ref{CGS-A:B} we also prove a weak converse to Theorem~\ref{CGS-Th:Bez_approx}.
Namely if $w.X_{\mathcal A}$ is sufficiently close to the geometric realization $|{\mathcal T}|$ of a
triangulation ${\mathcal T}$, then ${\mathcal T}$ is in fact the regular triangulation of ${\mathcal A}$
induced by the lifting function $\lambda({\bf a})=\log(w_{\bf a})$.

\begin{theorem}\label{CGS-Th:conv_Bez_approx}
  Let ${\mathcal A}\subset\R^d$ be a finite set of exponents.
  Suppose that $|{\mathcal T}|\subset \Sim^{\mathcal A}$ is the geometric realization of a triangulation
  ${\mathcal T}$ of ${\mathcal A}$ and there is a weight $w$ such that the distance between
  $X_{{\mathcal A},w}$ and $|{\mathcal T}|$ is less than $1/2(d+1)$. 
  Then ${\mathcal T}$ is the regular triangulation induced by the lifting function 
  $\lambda({\bf a})=\log(w_{\bf a})$.
\end{theorem}

\appendix
%
\section{Proofs of Injectivity Theorems}\label{CGS-A:A}
%

\noindent{\bf Theorem~\ref{CGS-Th:CF}.} (Craciun-Feinberg) {\it 
  The map $G_k$ is injective for every $k\in\R^m_>$ if and only 
  if the determinant of the Jacobian matrix,
\[
   \mbox{\rm Jac}(G_k)\ =\ 
    \left(\frac{\partial (G_k)_i}{\partial x_j}\right)_{i,j=1}^n\enspace ,
\]
  does not vanish for any $x\in\R^n_>$ and any $k\in\R^m_>$.
}\medskip

\begin{proof}
 We show the equivalence of the two statements, transforming one into the other.
 First, suppose there is a $k\in\R^m_>$ so that $G_k$ is not injective.
 Then there exist $c,d\in\R^n_>$ so that $G_k(c)=G_k(d)$.
 Then we have
 \[
   0\ =\ \sum_{i=1}^m k_i(c^{{\bf y}_i}-d^{{\bf y}_i}) {\bf z}_i\ =\ 
         \sum_{i=1}^m k'_i((\tfrac{c}{d})^{{\bf y}_i}-1) {\bf z}_i\enspace ,
 \]
 where $k'\in\R^m_>$ is defined by $\DeCo{k'_i}=k_i d^{{\bf y}_i}$, and 
 $\frac{c}{d}\in\R^n_>$ has $i$th coordinate $\frac{c_i}{d_i}$.
 In particular, $G_{k'}(\frac{c}{d})=G_{k'}(\iota)$, where $\DeCo{\iota}:=(1,\dotsc,1)$.
 Define ${\bf v}\in\R^n$ by 
\[
   v_i\ :=\ \log(c_i)-\log(d_i)\enspace,
\]
 so that $e^{{\bf y}_i\cdot {\bf v}}=(\frac{c}{d})^{{\bf y}_i}$, where ${\bf y}_i\cdot {\bf v}$
 is the Euclidean dot product, and we now have
 \begin{equation}\label{CGS-Eq:not_injective}
   0\ =\  \sum_{i=1}^m k'_i(e^{{\bf y}_i\cdot {\bf v}}-1) {\bf z}_i\enspace.
 \end{equation}

 Define the univariate function by $\DeCo{f}(t)=(e^t-1)/t$ for $t\neq 0$ and 
 set $f(0)=e$.
 Then $f$ is an increasing continuous bijection between $\R$ and $\R_>$.
 Define $k''\in\R^m_>$ by $\DeCo{k''_i}:=k'_i f({\bf y}_i\cdot{\bf v})$.
 Then $k'_i(e^{{\bf y}_i\cdot{\bf v}}-1)=k''_i({\bf y}_i\cdot{\bf v})$, and~\eqref{CGS-Eq:not_injective} becomes 
 \begin{equation}\label{CGS-Eq:magic_subs}
   0\ =\  \sum_{i=1}^m k''_i({\bf y}_i\cdot{\bf v}) {\bf z}_i\enspace.
 \end{equation}

 We claim that ${\bf v}$ lies in the kernel of the
 Jacobian matrix of $G_{k''}$ evaluated at the point $\iota$.
 Indeed, let $\DeCo{{\bf e}_j}:=(0,\dotsc,1,\dotsc,0)\in\R^n$ be the unit vector in the $j$th 
 direction.
 Then
\[
    \frac{\partial G_{k''}}{\partial x_j} (x)\ =\ 
    \sum_{i=1}^m k''_i x^{{\bf y}_i-{\bf e}_j} y_{i,j} {\bf z}_i\enspace,
\]
 where ${\bf y}_i=(y_{i,1},\dotsc,y_{i,n})$.
 Since $\iota^{{\bf y}_i-{\bf e}_j}=1$ and $\sum_j y_i v_j={\bf y}_i\cdot{\bf v}$, we see that 
\[
  \sum_{j=1}^n \frac{\partial G_{k''}}{\partial x_j} (\iota) \, v_j\ =\ 
  \sum_{j=1}^n \sum_{i=1}^m k''_i\iota^{{\bf y}_i-{\bf e}_j}y_{i,j} v_j {\bf z}_i\ =\ 
  \sum_{i=1}^m k''_i({\bf y}_i\cdot {\bf v}) {\bf z}_i\ =\ 0\enspace,
\]
 so that ${\bf v}$ lies in the kernel of the Jacobian matrix of $G_{k''}$ evaluated at 
 $\iota$, which implies that the Jacobian determinant of $G_{k''}$ vanishes at $\iota$.

 The theorem follows as these arguments are reversible.
\end{proof}

\noindent{\bf Corollary~\ref{CGS-C:matroid}.} {\it 
 The map $G_k$~$\eqref{CGS-Eq:exp_sum}$ is injective for all $k\in\R^m_>$ 
 if and only if $(Y_I\cdot Z_I)\cdot(Y_J\cdot Z_J)\geq 0$ for every 
 $I,J\in\binom{[m]}{n}$ and at least one product $Y_I\cdot Z_I$ is non-zero.
}\medskip


\begin{proof}
  Observe first that the Jacobian matrix $\mbox{Jac}(G_k)$ factors
  as the product of matrices $\delta^{-1} YDZ^T$, where $\delta$ is the diagonal matrix 
  with entries $(x_1,\dotsc,x_n)$, $D$ is the diagonal matrix with 
  entries $(k_1 x^{{\bf y}_1}, \dotsc, k_m x^{{\bf y}_m})$ and $Y$ and $Z$ are the matrices
  whose columns are the vectors ${\bf y}_i$ and ${\bf z}_i$, respectively.
  If we apply the Binet-Cauchy Theorem to this factorization, we see that
 \begin{equation}\label{CGS-Eq:Cauchy-Binet}  
   \det\mbox{Jac}(G_k)\ =\ x^{-\iota}\cdot \sum_{I\in\binom{[m]}{n}} \prod_{i\in I} k_i x^{{\bf y}_i}
    \cdot Y_I\cdot Z_I\enspace,
 \end{equation}
  where $\iota=(1,\dotsc,1)$.

  Suppose that $(Y_I\cdot Z_I)\cdot(Y_J\cdot Z_J)\geq 0$ for every $I,J\in\binom{[m]}{n}$,
  and at least one product $Y_I\cdot Z_I$ is non-zero.
  Then all terms in the sum~\eqref{CGS-Eq:Cauchy-Binet} have the same sign and not all are
  zero, and so the Jacobian does not vanish for any $x\in\R^n_>$ and $k\in\R^m_>$.
  Thus $G_k$ is injective for all $k\in\R^m_>$, by Theorem~\ref{CGS-Th:CF}.

  Suppose that there are two subsets $I,J\in\binom{[m]}{n}$ such that
  $Y_IZ_I>0$ and  $Y_JZ_J<0$.
  For $t\in\R_>$ and $K\in\binom{[m]}{n}$, define $k(K,t)\in\R^m_>$ by
\[
   k(K,t)_j\ :=\ \left\{\begin{array}{rcl}
     t&\ &\mbox{if }j\in K\\ 
     1&\ &\mbox{otherwise}\end{array}\right.
\]
 If we fix $x\in\R^d_>$, then the expansion~\eqref{CGS-Eq:Cauchy-Binet}, implies that 
 $\det\mbox{Jac}(G_{k(K,t)})(x)$ has the same sign as $Y_KZ_K$ when $t\gg 0$, at least
 when $Y_KZ_K\neq 0$. 
 
 We conclude that there is some $k,x$ such that $\det\mbox{Jac}(G_k)(x)>0$ and 
 some $k,x$ such that $\det\mbox{Jac}(G_k)(x)<0$,
 and therefore some $k,x$ such that $\det\mbox{Jac}(G_k)(x)=0$.
 This implies the corollary.  
\end{proof}

%
%
\section{Three Toric Theorems}\label{CGS-A:B}


\noindent{\bf Theorem~\ref{CGS-Th:X_calA}.} {\it
  Suppose that ${\mathcal A}\subset\R^d$ is a finite set of points with convex hull $\Delta$.
  Let $\beta=\{\beta_{\bf a}\mid {\bf a}\in{\mathcal A}\}$ be a collection of irrational toric B\'ezier
  functions for ${\mathcal A}$. 
  Then $\beta(\Delta)=X_{\mathcal A}$, the closure of the image of $\varphi_{\mathcal A}$.
}\medskip

\begin{proof}
 Let $\DeCo{\Delta^\circ}$ be the interior of $\Delta$, which we assume has $\ell$ facets
 and is given by the facet inequalities $0 \leq h_i(x), i=1,\dotsc,\ell$.
 Define two maps $\DeCo{H}\colon\Delta^\circ\to\R_>^\ell$ and 
 $\DeCo{\psi}\colon\R^\ell_>\to\Sim^{\mathcal A}$ by
 \begin{eqnarray*}
   \DeCo{H}\ \colon\ x&\longmapsto& (h_1(x),\dotsc,h_\ell(x))\,,\\
   \DeCo{\psi}\ \colon\ u&\longmapsto& 
      [u_1^{h_1({\bf a})}\dotsb u_\ell^{h_\ell({\bf a})}\ :\ {\bf a}\in{\mathcal A}]\,.
 \end{eqnarray*}
 Then the map $\beta\colon\Delta^\circ\to\Sim^{\mathcal A}$ (whose image is dense in 
 $X_{\mathcal A}$) is the composition of the maps $H$ and $\psi$.
\[
   \beta\ \colon\ \Delta^\circ\ \xrightarrow{\ H\ }\ 
    \R^\ell\  \xrightarrow{\ \psi\ }\ \Sim^{\mathcal A}\enspace.
\]

 Let us recall the definiton of the map 
$\DeCo{\varphi_{\mathcal A}}\colon\R_>^d\to\Sim^{\mathcal A}$,
\[
   \varphi_{\mathcal A}\ \colon\ 
     (x_1,\dotsc,x_d)\ \longmapsto\ [ x^{\bf a}\ :\ {\bf a}\in{\mathcal A}]\enspace.
\]
 The theorem follows once we show that the map $\psi$ factors through the map
 $\varphi_{\mathcal A}$.
 For this, define a new map $f_\Delta\colon\R_>^\ell\to\R^d_>$ by
\[
  f_\Delta\ \colon\ u\ \longmapsto\ t=(t_1,\dotsc,t_d)\qquad\mbox{where}\quad
    t_j\ :=\ \prod_{i=1}^\ell u_i^{{\bf v}_i\cdot{\bf e}_j}\enspace.
\]
 Then we claim that 
\[
   \psi\ \colon\  \R^\ell\ \xrightarrow{\ f_\Delta\ }\ \R^d_>
           \ \xrightarrow{\ \varphi_{\mathcal A}\ }\ \Sim^{\mathcal A}\enspace.
\]
 To see this, we compute the component of $\psi(u)$ corresponding to ${\bf a}\in{\mathcal A}$,
\[
   \prod_{i=1}^\ell u_i^{h_i({\bf a})}\ =\ 
   \prod_{i=1}^\ell u_i^{{\bf v}_i\cdot {\bf a} + c_i}\ =\ 
   \prod_{i=1}^\ell u_i^{c_i}\ \cdot\ \prod_{i=1}^\ell u_i^{{\bf v}_i\cdot {\bf a}}\ =\ 
   u^c \cdot  t^{\bf a}\enspace.
\]
 Thus $\psi(u)=u^c\varphi_{\mathcal A}(f_\Delta(u))$, as maps to $\R^{\mathcal A}_>$.
 The common factor $u^c$ does not affect the image in $\Sim^{\mathcal A}$, which shows that
 $\psi=\varphi_{\mathcal A} \circ f_\Delta$ and proves the theorem.  
\end{proof}

A consequence of this proof of Theorem~\ref{CGS-Th:X_calA} is the derivation of equations
which define the points of $X_{\mathcal A}$.
This derivation is similar to, but easier than, the development of toric ideals
in~\cite[Ch.~4]{CGS-GBCP}, as we have monomials with arbitrary real-number exponents. 

Suppose that we have a linear relation among the points of ${\mathcal A}$,
 \begin{equation}\label{CGS-Eq:linear_rel}
   \sum_{{\bf a}\in{\mathcal A}} \mu_{\bf a} \cdot {\bf a}\ =\ 
   \sum_{{\bf a}\in{\mathcal A}} \nu_{\bf a} \cdot {\bf a}\enspace,
 \end{equation}
for some $\mu,\nu\in\R^{\mathcal A}$.
Then the analytic binomial
 \begin{equation}\label{CGS-Eq:anal_binom}
   \prod_{{\bf a}\in{\mathcal A}} z_{\bf a}^{\mu_{\bf a}}\ -
   \prod_{{\bf a}\in{\mathcal A}} z_{\bf a}^{\nu_{\bf a}}
   \ =:\ \DeCo{z^\mu\ -\ z^\nu} 
 \end{equation}
vanishes on $\varphi_{\mathcal A}(\R^d_>)$, considered as a point in 
$\R^{\mathcal A}_>$.
This follows from the easy calculation
\[
  \varphi^*_{\mathcal A}(z^\mu)\ =\ 
  \prod_{{\bf a}\in{\mathcal A}} (x^{\bf a})^{\mu_{\bf a}}\ =\ 
  x^{\sum_{\bf a} \mu_{\bf a}\cdot {\bf a}}\enspace .
\]

Even after clearing denominators, the common zero set of the binomials~\eqref{CGS-Eq:anal_binom} 
is not exactly the image $\varphi_{\mathcal A}(\R^d_>)$ in the simplex $\Sim^{\mathcal A}$,
as the point $\varphi_{\mathcal A}(x)\in\Sim^{\mathcal A}$ is where the ray
\[
   \R_>\cdot (x^{\bf a}\mid {\bf a}\in{\mathcal A}) \ \subset\ \R^{\mathcal A}_>
\]
meets the simplex $\Sim^{\mathcal A}$.
If we require that the binomial~\eqref{CGS-Eq:anal_binom} is homogeneous in that 
$\sum_{{\bf a}\in{\mathcal A}} \mu_{\bf a}=\sum_{{\bf a}\in{\mathcal A}}\nu_{\bf a}$, then it
vanishes at every point of 
this ray and therefore on the image of $\varphi_{\mathcal A}$ in $\Sim^{\mathcal A}$.
Since the coordinates are positive numbers, we may further assume that~\eqref{CGS-Eq:linear_rel} is
an affine relation in that  
\[
   \sum_{{\bf a}\in{\mathcal A}} \mu_{\bf a}\ =\ \sum_{{\bf a}\in{\mathcal A}} \nu_{\bf a} \ =\ 1\enspace.
\]
These necessary conditions are also sufficient.

\begin{proposition}\label{CGS-Prop:eqs}
  A point $z$ in $\Sim^{\mathcal A}$ lies in $X_{\mathcal A}$ if and only if we have
\[
    z^\mu\ -\ z^\nu\ =\ 0
\]
  for all $\mu,\nu\in\R^{\mathcal A}$ with $\sum_{{\bf a}\in{\mathcal A}}\mu_{\bf a}=\sum_{{\bf a}\in{\mathcal A}}\nu_{\bf a}=1$ and
\[
   \sum_{{\bf a}\in{\mathcal A}} \mu_{\bf a}\cdot {\bf a}\ =\ 
   \sum_{{\bf a}\in{\mathcal A}} \nu_{\bf a}\cdot {\bf a}\enspace .
\]
\end{proposition}

One way to see the sufficiency is to pick an affinely independent subset ${\mathcal C}$ of
${\mathcal A}$ that affinely spans $\R^d$ and use the formula $x^{\bf a}=z_{\bf a}$ for ${\bf a}\in{\mathcal C}$ to
solve for $x\in\R^d$.
Then the point $z\in\Sim^{\mathcal A}$ satisfies this collection of binomials if and only if
$z=\varphi_{\mathcal A}(x)$.
This is also evident if we take logarithms of the coordinates.

These arguments only work for points $z$ in the interior of $\Sim^{\mathcal A}$.
For points of $X_{\mathcal A}$ with some coordinates zero, we use the recursive nature of
polytopes and the toric B\'ezier functions.
Namely, if we restrict the collection $\{\beta_{\bf a}\mid{\bf a}\in{\mathcal A}\}$ of toric B\'ezier
functions to a face ${\mathcal F}$ of the convex hull of ${\mathcal A}$, then those whose index ${\bf a}$
does not lie in ${\mathcal F}$ vanish, while those indexed by points of ${\mathcal A}$ lying in ${\mathcal F}$
specialize to toric B\'ezier functions for ${\mathcal F}$.
\medskip

\noindent{\bf Theorem~\ref{CGS-Th:Bez_approx}.} {\it
  Let ${\mathcal A}\subset\R^d$,  $w\in\R^{\mathcal A}_>$, and ${\mathcal B}\subset\R^n$  be exponents,
  weights, and control points for a toric B\'ezier patch.
  Suppose that ${\mathcal T}$ is a regular triangulation of ${\mathcal A}$ induced by a lifting function
  $\lambda\colon{\mathcal A}\to\R$. 
  For each $t>0$, let $F_t\colon \Delta\to\R^n$ be the toric B\'ezier patch of shape
  ${\mathcal A}$ with control points ${\mathcal B}$ and weights $t^{\lambda({\bf a})}w_{\bf a}$.
  Then, for any $\epsilon>0$ there exists a $t_0$ such that if $t>t_0$, the image 
  $F_t(\Delta)$ lies within $\epsilon$ of the control polytope ${\mathcal B}({\mathcal T})$.\medskip
}

\begin{proof}
 The lifting function $\lambda\colon{\mathcal A}\to\R$ inducing the triangulation ${\mathcal T}$ also 
 induces an action of $\R_>$ on $\Sim^{\mathcal A}$ where $t\in\R_>$ acts on a point
 $z\in\Sim^{\mathcal A}$ by scaling its coordinates, $(t.z)_{\bf a}=t^{\lambda({\bf a})} z_{\bf a}$.
 Then $F_t(\Delta)$ is the image of $t.X_{{\mathcal A},w}$ under the projection
 $\pi_{\mathcal B}\colon\Sim^{\mathcal A}\to\R^n$. 
 It suffices to show that $t.X_{{\mathcal A},w}$ can be made
 arbtrarily close to $|{\mathcal T}|\subset\Sim^{\mathcal A}$, if we choose $t$ large enough.

 We single out some equations from Proposition~\ref{CGS-Prop:eqs}.
 Suppose that ${\bf a},{\bf b}\in{\mathcal A}$ are points that do not lie in a common simplex of ${\mathcal T}$.
 That is, the segment $\overline{{\bf a},{\bf b}}$ is not an edge in the triangulation ${\mathcal T}$,
 and therefore it meets the interior of some face ${\mathcal F}$ of ${\mathcal T}$ so that 
 there is a point common to the interiors of ${\mathcal F}$ and of $\overline{{\bf a},{\bf b}}$.
 (If ${\bf a}={\bf b}$, so that ${\bf a}$ does not participate in the triangulation ${\mathcal T}$, then this
 point is just ${\bf a}$.)
 This gives the equality of convex combinations
 \begin{equation}\label{CGS-eq:convex_combination}
   \mu {\bf a} + (1{-}\mu){\bf b}\ =\ \sum_{{\bf c}\in{\mathcal F}} \nu_{\bf c}\cdot {\bf c}\enspace ,
 \end{equation}
 where all coefficients are positive and $\sum_{\bf c} \nu_{\bf c}=1$.
 Thus
\[
   z_{\bf a}^\mu z_{\bf b}^{1-\mu}\ =\ \prod_{{\bf c}\in{\mathcal F}} z_{\bf c}^{\nu_{\bf c}}
\]
 holds on $X_{\mathcal A}$.
 The corresponding equation on $t.X_{{\mathcal A},w}$ is
 \begin{equation}\label{CGS-eq:full}
     z_{\bf a}^\mu z_{\bf b}^{1-\mu}\ =\ 
     t^{\mu{\bf a} + (1{-}\mu){\bf b} - \sum_{{\bf c}\in{\mathcal F}} \nu_{\bf c}\cdot {\bf c}}
     \cdot \frac{w_{\bf a}^\mu w_{\bf b}^{1-\mu}}{\prod_{{\bf c}\in{\mathcal F}} w_{\bf c}^{\nu_{\bf c}}}
     \cdot   \prod_{{\bf c}\in{\mathcal F}} z_{\bf c}^{\nu_{\bf c}}\enspace .
 \end{equation}

Since $\overline{{\bf a},{\bf b}}$ is not in the triangulation, points in the
interior of the lifted segment
\[
   \overline{\,({\bf a},\lambda({\bf a})),\,({\bf b},\lambda({\bf b}))\,}
\]
lie below points of upper faces of the polytope $P_\lambda$.
We apply this observation to the point~\eqref{CGS-eq:convex_combination}.
Its height in the lifted segment is $\mu{\bf a} + (1{-}\mu){\bf b}$, while its height in the lift
of the face ${\mathcal F}$ is $\sum_{{\bf c}\in{\mathcal F}} \nu_{\bf c}\cdot {\bf c}$, and so
\[
    \mu{\bf a} + (1{-}\mu){\bf b}\ <\ \sum_{{\bf c}\in{\mathcal F}} \nu_{\bf c}\cdot {\bf c}\enspace .
\]
This implies that the exponent of $t$ in~\eqref{CGS-eq:full} is negative.
Since the other terms on the right hand side are bounded, we see that the left hand side,
and in fact the simple product $z_{\bf a} z_{\bf b}$, may be made as small as we please by
requiring that $t$ be sufficiently large.

Suppose that ${\mathcal A}$ consists of $\ell+d+1$ elements.
Repeating the previous argument for the (finitely many) pairs of points ${\bf a},{\bf b}$ which are not
both in any simplex of ${\mathcal T}$,
we see that for any $\epsilon>0$, there is a $t_0$ such that if $t>t_0$ and 
$z\in X_{{\mathcal A},w}$, then 
\[
   z_{\bf a} z_{\bf b} \ <\ \epsilon^2/4\ell^2\enspace,
\]
whenever ${\bf a},{\bf b}$ do not lie in a common simplex of ${\mathcal T}$.
In particular, at most one of $z_{\bf a}$ or $z_{\bf b}$ can exceed $\epsilon/2\ell$.

Let $z\in X_{{\mathcal A},w}$.  
Then there is some facet ${\mathcal F}$ of ${\mathcal T}$ such that if ${\bf a}\not\in{\mathcal F}$,
then $0\leq z_{\bf a}<\epsilon/2\ell$.
Suppose that ${\mathcal F}=\{{\bf a}_0,{\bf a}_1,\dotsc,{\bf a}_d\}$ and set
\[
   \DeCo{z_{\mathcal F}}\ :=\ (1-z_{{\bf a}_1}-\dotsb-z_{{\bf a}_d})e_{{\bf a}_0}
    + z_{{\bf a}_1}e_{{\bf a}_1} + \dotsb + z_{{\bf a}_d}e_{{\bf a}_d}\enspace ,
\]
which is a point of the facet $|{\mathcal F}|$ of the geometric realization $|{\mathcal T}|\subset\Sim^{\mathcal A}$.
Then
 \begin{eqnarray*}
   \| z-z_{\mathcal F}\| &=& 
   \left\| \sum_{{\bf a}\in{\mathcal A}} z_{\bf a} e_{\bf a}\ -\ 
     (1-z_{{\bf a}_1}-\dotsb-z_{{\bf a}_d})e_{{\bf a}_0}
    - z_{{\bf a}_1}e_{{\bf a}_1} - \dotsb - z_{{\bf a}_d}e_{{\bf a}_d}\right\|\\
   &\leq& \sum_{{\bf a}\not\in{\mathcal F}} z_{\bf a}\ +\ (1-z_{{\bf a}_0}-\dotsb-z_{{\bf a}_d})
    \ \ =\ 2 \sum_{{\bf a}\not\in{\mathcal F}} z_{\bf a} \\
   &<& 2\ell \frac{\epsilon}{2\ell}\ =\ \epsilon\,,
 \end{eqnarray*}
 as $1=\sum_{\bf a} z_{\bf a}$.
\end{proof}

We say that two subsets $X$ and $Y$ of Euclidean space are within a distance $\epsilon$ if for every
point $x$ of $X$ there is some point $y$ of $Y$ within a distance $\epsilon$ of $x$, and
vice-versa.\medskip 

\noindent{\bf Theorem~\ref{CGS-Th:conv_Bez_approx}. }{\it
  Let ${\mathcal A}\subset\R^d$ be a finite set of exponents.
  Suppose that $|{\mathcal T}|\subset \Sim^{\mathcal A}$ is the geometric realization of a triangulation
  ${\mathcal T}$ of ${\mathcal A}$ and there is a weight $w$ such that the distance between
  $X_{{\mathcal A},w}$ and $|{\mathcal T}|$ is less than $1/2(d+1)$. 
  Then ${\mathcal T}$ is the regular triangulation induced by the lifting function 
  $\lambda({\bf a})=\log(w_{\bf a})$.\medskip
}

\begin{proof}
 To show that ${\mathcal T}$ is the regular triangulation induced by the lifting function
 $\lambda$ whose value at ${\bf a}$ is $\log(w_{\bf a})$, we must show that if a segment
 $\overline{{\bf a},{\bf b}}$ between two points of ${\mathcal A}$ does not lie in the triangulation
 ${\mathcal T}$, then its lift by $\lambda$ lies below the lift of some face ${\mathcal F}$ of ${\mathcal T}$.

 Set $\epsilon:=1/2(d{+}1)$. 
 For each face ${\mathcal F}$ of ${\mathcal T}$, let $\DeCo{x_{\mathcal F}}\in\Sim^{\mathcal A}$ be the barycenter of
 ${\mathcal F}$,
\[
    x_{\mathcal F}\ :=\ \sum_{{\bf a}\in{\mathcal F}} \frac{1}{\#{\mathcal F}} e_{\bf a}\enspace ,
\]
 where $\#{\mathcal F}$ is the number of points of ${\mathcal A}$ in ${\mathcal F}$, which is at
 most $d{+}1$. 
 If $z$ is a point of $X_{{\mathcal A},w}$ within a distance $\epsilon$ of $x_{\mathcal F}$, 
 so that $\|z-x_{\mathcal F}\|<\epsilon$, then in particular no component of the vector
 $z-x_{\mathcal F}$ has absolute value exceeding $\epsilon$.
 Thus we have the dichotomy 
 \begin{equation}\label{CGS-Eq:dichotomy}
   \begin{array}{rrclcl}z_{\bf a}&<&\epsilon\ =\ 1/2(d{+}1)&\ &\mbox{if }\ {\bf a}\not\in{\mathcal F}\,,\\
                  z_{\bf a}&>&1/\#{\mathcal F}-\epsilon>1/2(d{+}1)&\ &\mbox{if }\ {\bf a}\in{\mathcal F}\enspace.
   \end{array}
\end{equation}

 Now suppose that the segment $\overline{{\bf a},{\bf b}}$ does not lie in the 
 triangulation ${\mathcal T}$.
 Then there is a face ${\mathcal F}$ of the triangulation whose interior meets the interior of
 this segment.
 That is, there is an equality of convex combinations~\eqref{CGS-eq:convex_combination} and a
 corresponding equation that holds for points $z\in X_{{\mathcal A},w}$,
\[
     z_{\bf a}^\mu z_{\bf b}^{1-\mu} \cdot \prod_{{\bf c}\in{\mathcal F}} w_{\bf c}^{\nu_{\bf c}}\ =\ 
     w_{\bf a}^\mu w_{\bf b}^{1-\mu} \cdot 
      \prod_{{\bf c}\in{\mathcal F}} z_{\bf c}^{\nu_{\bf c}}\enspace .
\]
 Suppose that $z$ is a point of $X_{{\mathcal A},w}$ that lies within a distance of $1/2(d{+}1)$
 of the barycenter $X_{\mathcal F}$ of the face ${\mathcal F}$.
 Then, by the estimates~\eqref{CGS-Eq:dichotomy}, we have
\[
  \frac{1}{2(d{+}1)} \prod_{{\bf c}\in{\mathcal F}} w_{\bf c}^{\nu_{\bf c}}\ >\ 
  z_{\bf a}^\mu z_{\bf b}^{1-\mu} \cdot \prod_{{\bf c}\in{\mathcal F}} w_{\bf c}^{\nu_{\bf c}}\ =\ 
     w_{\bf a}^\mu w_{\bf b}^{1-\mu} \cdot 
      \prod_{{\bf c}\in{\mathcal F}} z_{\bf c}^{\nu_{\bf c}}
  \ >\ w_{\bf a}^\mu w_{\bf b}^{1-\mu} \cdot \frac{1}{2(d{+}1)}\enspace .
\]
 canceling the common factor of $1/2(d{+}1)$ and taking logarithms, we obtain
\[
   \sum_{{\bf c}\in{\mathcal F}} \nu_{\bf c} \log(w_{\bf c})\ >\ 
   \mu \log(w_{\bf a}) + (1-\mu)\log(w_{\bf b})\enspace,
\]
 which implies that the point~\eqref{CGS-eq:convex_combination} common to the segment
 $\overline{{\bf a},{\bf b}}$ and the face ${\mathcal F}$ of ${\mathcal T}$ is lifted higher in the face
 ${\mathcal F}$ than in the segment $\overline{{\bf a},{\bf b}}$, and so the lift of the segment
 $\overline{{\bf a},{\bf b}}$ by $\lambda$ lies below the lift of the face ${\mathcal F}$.
 As this is true for all segments, we see that ${\mathcal T}$ is the triangulation induced by the 
 lifting function $\lambda$.
\end{proof}
\providecommand{\bysame}{\leavevmode\hbox to3em{\hrulefill}\thinspace}
\providecommand{\MR}{\relax\ifhmode\unskip\space\fi MR }
\providecommand{\MRhref}[2]{%
  \href{http://www.ams.org/mathscinet-getitem?mr=#1}{#2}
}
\providecommand{\href}[2]{#2}

\end{document}